\documentclass[pdflatex,sn-mathphys-num]{sn-jnl}
\usepackage{graphicx}
\usepackage{multirow}
\usepackage{amsmath,amssymb,amsfonts}
\usepackage{amsthm}
\usepackage[title]{appendix}
\usepackage{xcolor}
\usepackage{textcomp}
\usepackage{manyfoot}
\usepackage{booktabs}
\usepackage{algorithm}
\usepackage{algorithmicx}
\usepackage{algpseudocode}
\usepackage{listings}
\usepackage{orcidlink}
\usepackage{hyperref}
\usepackage{cleveref}
\crefname{figure}{Fig.}{Figs.}
\crefname{section}{Sec.\@}{Sec.\@}
\crefname{appendix}{Appendix}{Appendices}
\crefname{theorem}{Theorem}{Theorems}
\crefname{proposition}{Proposition}{Propositions}
\crefname{corollary}{Corollary}{Corollaries}
\crefname{remark}{Remark}{Remarks}

\newcommand{\R}{\mathbb{R}}
\newcommand{\cN}{\mathcal{N}}
\newcommand{\PrCov}{\mathbf\Gamma_{\mathrm{pr}}}
\newcommand{\ObsCov}{\mathbf\Gamma_{\mathrm{obs}}}
\newcommand{\PosCov}{\mathbf\Gamma_{\mathrm{pos}}}

\theoremstyle{thmstyleone}
\newtheorem{theorem}{Theorem}[section]
\newtheorem{proposition}[theorem]{Proposition}
\newtheorem{corollary}[theorem]{Corollary}
\newtheorem{remark}[theorem]{Remark}

\begin{document}
\title{Dimension and model reduction approaches for linear Bayesian inverse problems with rank-deficient prior covariances}

\author*[1]{\fnm{Josie} \sur{König}\orcidlink{0000-0003-4999-6649}} \email{josie.koenig@uni-potsdam.de}
\author[2]{\fnm{Elizabeth} \sur{Qian}\orcidlink{0000-0001-6713-3746}}\email{eqian@gatech.edu}
\author[1]{\fnm{Melina A.} \sur{Freitag}\orcidlink{0000-0002-4539-2162}} \email{melina.freitag@uni-potsdam.de}
\affil*[1]{ \orgname{Universität Potsdam}, \orgdiv{Institut für Mathematik}, \orgaddress{\street{Karl-Liebknecht-Str. 24--25}, \postcode{14476} \city{Potsdam}, \country{Germany}}}
\affil[2]{\orgname{Georgia Institute of Technology}, \orgaddress{\street{North Avenue}, \city{Atlanta},  \state{GA} \postcode{30332}, \country{USA}}}

\abstract{
Bayesian inverse problems use observed data to update a prior probability distribution for an unknown state or parameter of a scientific system to a posterior distribution conditioned on the data.
In many applications, the unknown parameter is high-dimensional, making computation of the posterior expensive due to the need to sample in a high-dimensional space and the need to evaluate an expensive high-dimensional forward model relating the unknown parameter to the data. 
However, inverse problems often exhibit low-dimensional structure due to the fact that the available data are only informative in a low-dimensional subspace of the parameter space. 
Dimension reduction approaches exploit this structure by restricting inference to the low-dimensional subspace informed by the data, which can be sampled more efficiently.
Further computational cost reductions can be achieved by replacing expensive high-dimensional forward models with cheaper lower-dimensional reduced models.
In this work, we propose new dimension and model reduction approaches for linear Bayesian inverse problems with rank-deficient prior covariances, which arise in many practical inference settings. The dimension reduction approach is applicable to general linear Bayesian inverse problems whereas the model reduction approaches are specific to the problem of inferring the initial condition of a linear dynamical system. We provide theoretical approximation guarantees as well as numerical experiments demonstrating the accuracy and efficiency of the proposed approaches.

}


\maketitle

\section{Introduction}
Inverse problems seek to infer an unknown state or parameter of a scientific system from observed data, and arise in many scientific disciplines, including meteorology \cite{Buehner2017VarEnKF,Houtekamer2005AtmoDAEnKF,Navon2009DAforNWP}, oceanography, \cite{BRASSEUR1999SEEK,Rozier2007ROKF,TUANPHAM1998SEEK,Verlaan1997RRSRKF} and environmental science \cite{Carrassi2018DAGeo,KITANIDIS2015CSEK,Leeuwen2010NonLinDAGeo}. Inverse problems are typically defined by an observation model which relates the observed data $\mathbf{m} \in \mathbb{R}^{d_\mathrm{obs}}$ to an unknown state or parameter $\mathbf{p}\in\mathbb{R}^d$ through the forward model $\mathbf{G}:\mathbb{R}^d \rightarrow \mathbb{R}^{d_\mathrm{obs}}$ that is subject to additive noise $\pmb{\epsilon}\in \mathbb{R}^{d_\mathrm{obs}}$:
$$\mathbf{m} = \mathbf{G}(\mathbf{p})+\pmb{\epsilon}.$$
In the Bayesian approach to inference \cite{Dashti2017BIPSurvey,Law2015DA,Taran2005BIPbook}, the prior distribution on $\mathbf{p}$, describing assumptions and associated uncertainties before any data are observed, is updated to the posterior distribution $\mathbf{p}|\mathbf{m}$ of the unknown parameter conditioned on the data using Bayes' rule.
In large scale scientific applications, computing this posterior distribution is often expensive due to the high dimension of $\mathbf{p}$ and the expense of evaluating $\mathbf{G}$, which may need to be evaluated many times to sample from the posterior distribution, e.g., via Markov chain Monte Carlo methods \cite{Bardsley2015CUQfIP,heard2021}, or to estimate the parameter for different realizations of the data. 

\textit{Dimension} reduction approaches seek to alleviate this computational burden by exploiting low-rank structure in inverse problems, which arises from the fact that the available data are often primarily informative (relative to the prior) in a low-dimensional subspace. This may be due to the limited amount of observed data and/or smoothing properties of the forward model. 
Several works seek to exploit this structure for computational cost reduction by updating the posterior only in this low-dimensional \textit{likelihood-informed subspace} (LIS) ~\cite{flathFastAlgorithmsBayesian2011,bui-thanhComputationalFrameworkInfiniteDimensional2013,martinStochasticNewtonMCMC2012,Spantini2015Optimal}. In particular, the work in~\cite{Spantini2015Optimal} shows that updates in the directions of the leading generalized eigenvectors of the Hessian of the negative log-likelihood and the prior precision matrix yield optimal posterior approximations for linear Gaussian Bayesian inverse problems. This work has been extended to nonlinear \cite{Cui2014LISNonlinear} and non-Gaussian \cite{Cui2021DatafreeLIS,Zahm2022NonlinCertifiedMOR} Bayesian inference problems and to goal-oriented inference problems \cite{Lieberman2013IfP}. These dimension reduction approaches yield computational cost reductions because the low-dimensional LIS can be more efficiently explored than the original high-dimensional unknown space, but dimension reduction approaches alone still require evaluation of the original high-dimensional forward model to sample in the low-dimensional LIS.

In contrast, \textit{model} reduction methods seek computational cost reductions by replacing the original high-dimensional forward model with a cheaper low-dimensional reduced model. Projection-based reduced models \cite{Antoulas2005Book,Benner2015Survey,HestRozStamm2015RBM-PDE} obtain a reduced model by projecting the governing equations of the high-dimensional system onto a low-dimensional subspace. In the popular proper orthogonal decomposition (POD) method, this reduced subspace is defined as the span of the leading principal components of available state data~\cite{Galbally2010NonLinMORPODInterpol,Lieberman2010GreedyPODInverseProb,Lipponen2013PODtomography,Nguyen2014PODBayesHydro,Wang2005POD-BIP-heat}. 
In contrast, system- or control-theoretic methods for model reduction define the reduced subspace not from data but from the structure of the system: such methods include rational interpolation methods~\cite{AntBeaGuc2020InterpolatoryMeths,Benner2017MORBook,GugercinIRKA2008}, which define the projection subspace so that the transfer function of the reduced system interpolates that of the full system, as well as balanced truncation~\cite{Moore1981BT,Mullis1976BT}, which projects the system onto state directions that are the leading generalized eigenvectors of the system observability and inverse reachability Gramians, preserving directions that are simultaneously highly observable and easily reachable. 
Recently, several works have developed \textit{inference-oriented} reduced models which exploit connections between system-theoretic balanced truncation methods and LIS dimension reduction: the work~\cite{Qian2021Balancing} showed that system-theoretic balanced truncation can be adapted to the inference problem of inferring the initial condition of a linear dynamical system from noisy measurements taken after the initial time, yielding an inference-oriented reduced model that recovers the optimal posterior approximation guarantees of~\cite{Spantini2015Optimal} for stable systems in certain limits. This approach has been generalized to unstable systems in~\cite{Koe2022TLBTDA,koenig2022TLBT4DVAR} and quadratic dynamical systems in~\cite{freitag2024Q-BT}.

To date, both LIS dimension reduction approaches as well as the related inference-oriented model reduction approaches have assumed that the prior covariance matrix is invertible in both the definition of the approximation and the development of associated optimality guarantees. 
However, rank-deficient prior covariances arise in many practical settings, including in particle filtering algorithms where the rank of the covariance is restricted by the size of the ensemble~\cite{Buehner2017VarEnKF,Houtekamer2005AtmoDAEnKF,LIVINGS2008UnbiasedEnKF,Leeuwen2010NonLinDAGeo}, as well as in large-scale applications such as weather forecasting where prior covariances may be poorly conditioned and very close to rank-deficient \cite{Buehner2010ErrorStatsDA,HABEN2011PreCond4dVar} or even approximated by low-rank factors \cite{BRASSEUR1999SEEK,KITANIDIS2015CSEK,Rozier2007ROKF,TUANPHAM1998SEEK,Verlaan1997RRSRKF}. 
The purpose of this work is to define dimension and model reduction approaches and theory for linear Bayesian inverse problems with rank-deficient (non-invertible) prior covariance matrices. Our contributions are:
\begin{enumerate}
    \item For both the LIS dimension reduction approach of~\cite{Spantini2015Optimal} as well as the LIS-balanced truncation model reduction approach of~\cite{Qian2021Balancing}, we provide new definitions of the posterior approximations for rank-deficient prior covariances that recover the approximations of~\cite{Spantini2015Optimal,Qian2021Balancing} if the prior covariance is invertible. We show that these approximations inherit optimality guarantees in metrics defined on a lower-dimensional space associated with the range of the prior covariance.
    \item We propose a new inference-oriented model reduction approach applicable to the inference of the initial condition of a linear dynamical system from linear measurements taken after the initial time. The new approach, \textit{prior-driven balancing} (PD-BT), is based on interpreting the unknown initial condition as an impulse input to a linear time-invariant control system whose input port matrix is given by a square root factor of the prior covariance matrix. We show that the PD-BT approach inherits system-theoretic stability guarantees and error bounds.
    \item We provide numerical experiments that demonstrate the accuracy and efficiency of the proposed methods and illustrate our theoretical results.
\end{enumerate}
The remainder of this paper is organized as follows: \Cref{sec: Background} introduces the problem formulation and relevant background. \Cref{sec: whole_sec_rankdef} introduces our new approximation definitions for rank-deficient prior covariance matrices. \Cref{sec: Numerics} contains numerical experiments and discussions and \Cref{sec: Conclusions} concludes.

\section{Background}\label{sec: Background}
We now introduce the problem formulation in \cref{subsec: problem formulation} before providing background on LIS dimension reduction in \cref{subsec: Spantini} and on balanced truncation methods for model reduction in \cref{subsec: BT}.
\subsection{Problem setting}\label{subsec: problem formulation}
We consider a Bayesian inverse problem with linear forward operator $\mathbf{G}\in \mathbb{R}^{{d_\mathrm{obs}}\times d}$ and centered Gaussian prior and noise distributions:
\begin{align}\label{eq: generic_IVP}
    \mathbf{m} = \mathbf{G}\mathbf{p}+\pmb{\epsilon}, \enspace\mathbf{p}\sim \mathcal{N}(\mathbf{0},\PrCov),\enspace\pmb{\epsilon}\sim \mathcal{N}(\mathbf{0},\ObsCov).
\end{align}
The posterior distribution for this setting is Gaussian \cite{Carlin2009BayesMeth}: $\mathbf{p}|\mathbf{m}\sim\mathcal{N}(\pmb{\mu}_{\mathrm{pos}},\PosCov)$ with
\begin{align}\label{eq: Posterior}
\begin{split}
    \pmb{\mu}_{\mathrm{pos}} &= \PosCov\mathbf{G}^\top\ObsCov^{-1}\mathbf{m} \in \mathbb{R}^d,\quad \text{and} \\\PosCov &= \PrCov - \PrCov \mathbf{G}^\top(\ObsCov+\mathbf{G}\PrCov \mathbf{G}^\top)^{-1}\mathbf{G}\PrCov\in \mathbb{R}^{d \times d}.
\end{split}
\end{align}
\Cref{eq: generic_IVP,eq: Posterior} describe a general linear Gaussian Bayesian inverse problem considered by LIS dimension reduction \cite{Spantini2015Optimal} and which we consider in \cref{subsec: Spantini,subsec: Spantini-rankdef}. 

We now introduce a special case of the linear Gaussian Bayesian inverse problem in which the unknown parameter is the initial condition of a high-dimensional linear dynamical system, and the data arise from noisy measurements of a linear system output taken after the initial time. Let $\mathbf{A} \in \mathbb{R}^{d \times d}$ and $\mathbf{C} \in \mathbb{R}^{d_{\mathrm{out}} \times d}$ and consider the dynamical system with state $\mathbf{x}(t)\in\mathbb{R}^d$ and output $\mathbf{y}(t)\in\mathbb{R}^{d_{\rm out}}$:
\begin{subequations}\label{eq: BayesianSystem}
\begin{align}
\begin{split}\label{eq: BayesianSystem_nonoise}
    \dot{\mathbf{x}}(t)  &= \mathbf{A}\mathbf{x}(t),\enspace\quad \mathbf{x}(0)=\mathbf{p},\\
    \mathbf{y}(t) &= \mathbf{C}\mathbf{x}(t).
\end{split}
\end{align}
The data arise from $n$ output measurements $\mathbf{m} = \begin{bmatrix} \mathbf{m}_1^\top, \ldots, \mathbf{m}_n^\top \end{bmatrix}^\top\in \mathbb{R}^{nd_{\mathrm{out}}}$ taken at positive times $t_{1} < \ldots < t_{n}$, subject to additive Gaussian noise with mean zero and positive definite covariance $\mathbf{\Gamma}_{\epsilon}\in\mathbb{R}^{d_{\rm out}\times d_{\rm out}}$:
\begin{align}\label{eq: BayesianSystem_noise}
    \mathbf{m}_k &= \mathbf{y}(t_k) + \epsilon_k= \mathbf{C}\mathbf{x}(t_k) + \epsilon_k, \enspace \epsilon_k \sim \mathcal{N}(\mathbf{0},\mathbf{\Gamma}_\epsilon), \enspace k=1,\ldots,n. 
\end{align}
\end{subequations}
In this work, we refer to the problem~\eqref{eq: BayesianSystem} of inferring the unknown initial condition $\mathbf{p}$ from measurements taken after $t=0$ as the \textit{Bayesian smoothing problem}, as commonly used in the Bayesian inverse problems literature \cite{Gelb1974, Law2015DA}.
Note that the Bayesian smoothing problem has the form~\eqref{eq: generic_IVP} and the solution~\eqref{eq: Posterior} with 
\begin{align*}
\mathbf{G}= \begin{bmatrix}\mathbf{C}e^{\mathbf{A}t_1}\\\vdots\\\mathbf{C}e^{\mathbf{A}t_n}\end{bmatrix}\in \mathbb{R}^{d_{\rm obs}\times d}, \quad\text{and}\quad \ObsCov= \begin{bmatrix}\mathbf{\Gamma}_{\pmb{\epsilon}}\\ & \ddots \\ &&\mathbf{\Gamma}_{\pmb{\epsilon}} \end{bmatrix}\in \mathbb{R}^{d_{\rm obs}\times d_{\rm obs}},
\end{align*}
where $d_{\rm obs} = n\cdot d_{\rm out}$.
In large-scale applications, it is often impractical to assemble the operator $\mathbf{G}$; instead the mapping between the initial condition and the measurements is only implicitly available through the time integration of the high-dimensional dynamical system~\eqref{eq: BayesianSystem_nonoise}. To reduce the cost of evaluating $\mathbf{G}$ in such settings, we will consider \textit{model} reduction approaches which replace the high-dimensional model~\eqref{eq: BayesianSystem_nonoise} with a lower-dimensional \textit{reduced model} in \Cref{subsec: LIS-BT,sec: Qian_rankdef,sec: PD-BT}.

\subsection{Likelihood-informed dimension reduction for Bayesian inverse problems}\label{subsec: Spantini}
This section introduces the likelihood-informed dimension reduction approach defined and analyzed in~\cite{flathFastAlgorithmsBayesian2011,Spantini2015Optimal} for general linear Gaussian Bayesian inverse problems of the form~\eqref{eq: generic_IVP}. The main idea is to preserve directions of the unknown parameter space which maximize the following Rayleigh quotient:
\begin{align*}
    \frac{\mathbf{p}^\top\mathbf (\mathbf{G}^\top\ObsCov^{-1}\mathbf G)\mathbf{p}}{\mathbf{p}^\top\PrCov^{-1}\mathbf{p}}.
\end{align*}
The quantity $\mathbf G^\top\ObsCov^{-1}\mathbf G$ is called the Fisher information matrix, and is the Hessian of the negative log-likelihood. Directions which maximize this quotient are highly informed by the data relative to their prior uncertainty, and are given by the leading generalized eigenvectors of the matrix pencil $(\mathbf G^\top\ObsCov^{-1}\mathbf G,\PrCov^{-1})$. Let: 
\begin{align}\label{eq: spantini GEV}
        \mathbf G^\top\ObsCov^{-1}\mathbf G\mathbf{W}=\PrCov^{-1}\mathbf{W}\mathbf{\Sigma}^2
    \end{align}
where $\mathbf{\Sigma}^2=\mathrm{diag}(\mathbf{\sigma}_1^2,\ldots,\mathbf{\sigma}_d^2)$ is the matrix of eigenvalues with $\mathbf{\sigma}_1^2\geq\mathbf{\sigma}_2^2\geq\ldots\geq\mathbf{\sigma}_d^2$ in non-increasing order. Note that some works, e.g.,~\cite{flathFastAlgorithmsBayesian2011}, consider an eigenvalue decomposition of the `prior-preconditioned Hessian of the data misfit', $\PrCov^{1/2}\mathbf G^\top\ObsCov^{-1}\mathbf G\PrCov^{1/2}$, whose eigenvalues are equal to the generalized eigenvalues of \eqref{eq: spantini GEV} and whose eigenvectors are related to the generalized eigenvectors of \eqref{eq: spantini GEV} via a linear transformation.

The matrix $\mathbf{W}$ collects the corresponding eigenvectors and without loss of generality, we impose the following normalization\footnote{This is a different normalization from the one presented in~\cite{Spantini2015Optimal}; we use this choice because it is consistent with the normalization most commonly used in balanced truncation model reduction.}: $\mathbf{W}^\top\mathbf G^\top\ObsCov^{-1}\mathbf G\mathbf{W}=\mathbf{\Sigma}$. Then, the matrix $\mathbf V :=\mathbf G^\top\ObsCov^{-1}\mathbf G\mathbf{W}\mathbf{\Sigma}^{-1}$ is orthogonal to $\mathbf{W}$, i.e., $\mathbf{V}^\top\mathbf{W}=\mathbf{I}_d$, and satisfies $\mathbf{V}^\top\PrCov\mathbf{V}=\mathbf{\Sigma}$.

Let $\mathbf{V}_r,\mathbf{W}_r \in\R^{d\times r}$ contain the first $r$ columns of $\mathbf V,\mathbf{W}$. Then, $\mathbf{\Pi}_r := \mathbf{W}_r \mathbf{V}_r^\top\in\R^{d\times d}$ is an oblique projector onto the $r$-dimensional \textit{likelihood-informed subspace} (LIS) spanned by the dominant generalized eigenvectors of~\eqref{eq: spantini GEV}. Projecting the unknown parameter $\mathbf{p}$ onto the LIS yields the following approximate measurement model: 
\begin{align*}
    \mathbf{m}=\mathbf{G}\mathbf{\Pi}_r\mathbf{p}+\pmb{\epsilon} \equiv \mathbf{G}_{\mathrm{OLR}}\,\mathbf{p}+\pmb{\epsilon},
\end{align*}
where $\mathbf{G}_{\mathrm{OLR}} := \mathbf G\mathbf{\Pi}_r \in \mathbb{R}^{d_{\rm obs}\times d}$ is a high-dimensional operator of rank (at most) $r$. This approximate measurement model then induces the following optimal low-rank (OLR) approximations to the posterior mean and covariance:
\begin{align}\label{eq: Spantini PosBoth}
\begin{split}
    \PosCov^{\mathrm{OLR}} &= (\PrCov^{-1} + \mathbf{G}_{\mathrm{OLR}}^\top\ObsCov^{-1}\mathbf{G}_{\mathrm{OLR}})^{-1}\\ &=  \PrCov - \PrCov \mathbf{G}_{\mathrm{OLR}}^\top(\ObsCov+\mathbf{G}_{\mathrm{OLR}}\PrCov \mathbf{G}_{\mathrm{OLR}}^\top)^{-1}\mathbf{G}_{\mathrm{OLR}}\PrCov,\\
    \pmb \mu_{\mathrm{pos}}^{\mathrm{OLR}} &=\PosCov^{\mathrm{OLR}}\mathbf{G}_{\mathrm{OLR}}^\top\ObsCov^{-1}\mathbf m.    
\end{split}
\end{align}
The work~\cite{Spantini2015Optimal} shows that the approximations~\eqref{eq: Spantini PosBoth} are optimal in the following sense:
\begin{theorem}[{\cite{Spantini2015Optimal}}]\label{thm: Spantini}
Let $\mathcal{M}_r$ denote the class of positive definite matrices in $\R^{d\times d}$ that are given by (at most) rank-$r$ negative semi-definite updates to the prior covariance, i.e., 
\begin{align*}
    \mathcal{M}_r =\{ \PrCov - \mathbf K_r\mathbf K_r^\top\succ 0: \mathrm{rank}(\mathbf K_r)\leq r\},
\end{align*}
and let $\mathcal{A}_r$ denote the class of vectors in $\R^d$ given by (at most) rank-$r$ operators acting on the data:
\begin{align*}
    \mathcal{A}_r= \{\mathbf A_r \mathbf{m} \in \mathbb{R}^d :\mathrm{rank}(\mathbf A_r)\leq r\}.
\end{align*}
Then, the OLR mean and covariance approximations~\eqref{eq: Spantini PosBoth} satisfy
\begin{align*}
    \PosCov^{\mathrm{OLR}}=\arg\min_{\hat{\mathbf{\Gamma}}\in\mathcal{M}_r}d^2_\mathcal{F}(\PosCov,\hat{\mathbf{\Gamma}}) && and && {\pmb \mu}_{\mathrm{pos}}^{\mathrm{OLR}} =\arg\min_{\hat{\pmb\mu}\in\mathcal{A}_r} \mathbb{E} \|\pmb\mu_{\rm pos} - \hat{\pmb\mu} \|_{\PosCov^{-1}}^2,
\end{align*}
where $d_{\mathcal{F}}$ denotes the F\"orstner distance between two symmetric positive-definite matrices~\cite{Foerstner2003}, defined as $d_{\mathcal{F}}^2 (\mathbf{E},\mathbf{F}):= \sum_i \ln^2(\lambda_i)$ where $\lambda_i$ are the generalized eigenvalues of the matrix pencil $(\mathbf{E},\mathbf{F})$, and the quantity $\mathbb{E} \|\pmb\mu_{\rm pos} - \hat{\pmb\mu} \|_{\PosCov^{-1}}^2$ is the Bayes risk induced by the Mahalanobis distance associated with the posterior distribution.
\end{theorem}

The optimal approximations~\eqref{eq: Spantini PosBoth} are more efficient to compute than the full posterior statistics~\eqref{eq: Posterior} in settings where the low-rank operator $\mathbf{G}_{\mathrm{OLR}}$ can be efficiently evaluated. However, computing $\mathbf{G}_{\mathrm{OLR}}$ may be impractical in settings where $\mathbf{G}$ involves the evolution of a large-scale dynamical system. The next section introduces methods for model reduction which yield cheap approximations for the dynamical system.

\subsection{Model order reduction by balanced truncation}\label{subsec: BT}
We now describe balanced truncation (BT) methods for model reduction, beginning with a general formulation (\Cref{subsubsec: general BT}). We then consider the canonical control-theoretic balanced truncation formulation for linear time-invariant (LTI) systems in \Cref{subsubsec: control-BT}, followed by an adaptation to the Bayesian smoothing problem for linear dynamical systems in \Cref{subsec: LIS-BT}. 

\subsubsection{General BT reduced model definition}\label{subsubsec: general BT}
Balanced truncation typically considers a high-dimensional linear time-invariant (LTI) system defined by matrix operators $\mathbf{A} \in \mathbb{R}^{d \times d}$, $\mathbf{B} \in \mathbb{R}^{d \times d_\mathrm{in}}$, $\mathbf{C} \in \mathbb{R}^{d_\mathrm{out} \times d}$ as follows:
    \begin{align}\label{eq: LTISys}
    \begin{split}
        \dot{\mathbf{x}}(t) &= \mathbf{Ax}(t) + \mathbf{Bu}(t),\\
        \mathbf{y}(t)  &= \mathbf{Cx}(t).
    \end{split}
\end{align}
For now, let $\mathbf{P},\mathbf{Q}\in\mathbb{R}^{d\times d}$ be symmetric positive definite matrices (we will discuss possible choices for these matrices shortly). Balanced truncation methods obtain a reduced model for~\eqref{eq: LTISys} via Petrov-Galerkin projection onto a low-dimensional subspace defined by $\mathbf{P}$ and $\mathbf{Q}$ in the following way: consider the generalized eigenvalue problem,
\begin{align}\label{eq: BT-EVP}
    \mathbf{Q}\mathbf{W}=\mathbf{P}^{-1}\mathbf{W}\mathbf{\Sigma}^2,
\end{align}
where $\mathbf{\Sigma}=\mathrm{diag}(\sigma_1,\sigma_2,\ldots,\sigma_d)\in\R^{d\times d}$ and $\mathbf{W} = \begin{bmatrix}\mathbf{w}_1 & \mathbf{w}_2 & \cdots & \mathbf{w}_d\end{bmatrix}\in\R^{d\times d}$ collect the $d$ pairs of generalized eigenvalues and eigenvectors of the matrix pencil $(\mathbf{Q},\mathbf{P}^{-1})$, with the eigenvalues ordered in non-increasing order $\sigma_1\geq \sigma_2 \geq \cdots \geq \sigma_d>0$ and the eigenvectors normalized so that $\mathbf{W}^\top\mathbf{Q}\mathbf{W} = \mathbf{\Sigma}$. Note that $\mathbf{V} = \mathbf{Q}\mathbf{W}\mathbf{\Sigma}^{-1}\in\R^{d\times d}$ satisfies $\mathbf{V}^\top\mathbf{W}=\mathbf{I}_d$: the column vectors of $\mathbf{V}$ are then the generalized eigenvectors of the matrix pencil $(\mathbf{P},\mathbf{Q}^{-1})$ normalized such that $\mathbf{V}^\top\mathbf{P}\mathbf{V} = \mathbf{\Sigma}$.
Let $\mathbf{W}_r,\mathbf{V}_r\in\R^{d\times r}$ denote matrices containing the leading $r$ columns of $\mathbf{W},\mathbf{V}$, respectively. These leading basis vectors maximize the Rayleigh quotients $\frac{\mathbf{x}^\top\mathbf{Q}\mathbf{x}}{\mathbf{x}^\top\mathbf{P}^{-1}\mathbf{x}}$ and $\frac{\mathbf{x}^\top\mathbf{P}\mathbf{x}}{\mathbf{x}^\top\mathbf{Q}^{-1}\mathbf{x}}$, respectively.

In balanced truncation, the state $\mathbf{x}$ is approximated in the span of the trial basis $\mathbf{W}_r$, that is, $\mathbf{x}\approx\mathbf{W}_r\mathbf{x}_r$, with $\mathbf{x}_r\in\R^r$ denoting the \textit{reduced} state. A reduced LTI system is then obtained by enforcing the Petrov-Galerkin orthogonality condition that the residual of the system dynamics must be orthogonal to the test basis $\mathbf{V}_r$, yielding 
\begin{align}\label{eq: RedLTISys}
    \begin{split}
        \dot{\mathbf{x}}_r(t) &= \mathbf{A}_r \mathbf{x}_r(t) + \mathbf{B}_r\mathbf{u}(t),\\
         \mathbf{y}_r(t)  &= \mathbf{C}_r \mathbf{x}_r(t),
    \end{split}
\end{align}
where $\mathbf{A}_r :=  \mathbf{V}_r^\top\mathbf{A}\mathbf{W}_r$, $\mathbf{B}_r :=  \mathbf{V}_r^\top \mathbf{B}$, and $\mathbf{C}_r := \mathbf{C}\mathbf{W}_r$.
The system~\eqref{eq: RedLTISys} is an $r$-dimensional representation of the projection of the $d$-dimensional system~\eqref{eq: LTISys} onto the range of $\mathbf{W}_r$ using the oblique projector $\mathbf{\Pi}_r = \mathbf{W}_r\mathbf{V}_r^\top$. This $r$-dimensional system can then be cheaply evolved to obtain approximations of the output signal $\mathbf{y}_r(t)$.

\begin{remark}\label{rem: square root balancing}
    Computing the balanced truncation reduced model by directly solving the generalized eigenvalue problem may suffer from ill-conditioning, especially if $\mathbf{P}$ and $\mathbf{Q}$ are near low-rank. A more numerically stable computational procedure is based on (low-rank) square root factors (e.g., Cholesky) of $\mathbf{P}$ and $\mathbf{Q}$. Let~$\mathbf{P} \approx \mathbf{LL}^\top, \mathbf{Q} \approx \mathbf{RR}^\top$, with $\mathbf{L},\mathbf{R}\in\R^{d\times s}$ with $s<d$. Then, for $r\leq s$, let $\mathbf{R}^\top \mathbf{L}\approx\mathbf{U}_{r}\mathbf{\Sigma}_{r}\mathbf{Z}_{r}^\top$ denote the rank-$r$ truncated singular value decomposition of the product $\mathbf{R}^\top \mathbf{L}$. Then, the left and right bases are given by $\mathbf{W}_r = \mathbf{L}\mathbf{Z}_{r}\mathbf{\Sigma}_{r}^{-\frac{1}{2}}$ and $\mathbf{V}_r =\mathbf{R}\mathbf{U}_{r}\mathbf{\Sigma}_{r}^{-\frac{1}{2}}$. Computing the balanced truncation bases in this way is called square-root balancing~\cite[p.\@ 221 ff.\@]{Antoulas2005Book}. An exactly analogous procedure for computing the likelihood-informed subspace bases is suggested in~\cite{Spantini2015Optimal}. In \Cref{sec: existing methods rankdef}, we discuss the interpretation of the basis definition via square-root balancing when $\mathbf{P}$ is a rank-deficient prior covariance, including the senses in which the optimality guarantees of \Cref{subsec: Spantini} are preserved.
\end{remark}

\subsubsection{Canonical balanced truncation for LTI control systems}\label{subsubsec: control-BT}
The canonical balanced truncation method for LTI control systems of the form~\eqref{eq: LTISys} assumes that the matrix $\mathbf{A}$ is stable, with eigenvalues in the open left-half plane, and is based on the following definitions of $\mathbf{P}$ and $\mathbf{Q}$:
\begin{align}\label{eq: Gramians}
    \mathbf{P} = \int_0^\infty e^{\mathbf{A}t}\mathbf{BB}^\top e^{\mathbf{A}^\top t}\,\text{d}t, \qquad \mathbf{Q} = \int_0^\infty e^{\mathbf{A}^\top t}\mathbf{C}^\top \mathbf{C} e^{\mathbf{A} t}\,\text{d}t,
\end{align}
which are called the (infinite) \textit{reachability} and \textit{observability} Gramians, respectively. The infinite Gramians are the unique solutions to the dual Lyapunov equations, $\mathbf{AP} + \mathbf{PA}^\top = -\mathbf{BB}^\top$ and $\mathbf{A}^\top \mathbf{Q} + \mathbf{QA} = -\mathbf{C}^\top\mathbf{C}$,
which can be efficiently solved for small- to medium-scale systems. For larger systems, low-rank Lyapunov solvers can be used to compute low-rank Cholesky factors of the Gramians~\cite{Benner2013NumSolSurvey} which may then be used to compute the BT bases using the square-root balancing procedure (see \Cref{rem: square root balancing}).

Recall that the state approximation basis $\mathbf{W}_r$ consists of the directions which maximize $\frac{\mathbf{x}^\top\mathbf{Q}\mathbf{x}}{\mathbf{x}^\top\mathbf{P}^{-1}\mathbf{x}}$. For the canonical system Gramians~\eqref{eq: Gramians}, the quadratic forms $\mathbf{x}^\top\mathbf{P}^{-1}\mathbf{x}$ and $\mathbf{x}^\top \mathbf{Q}\mathbf{x}$ define system reachability and observability energies, respectively. The reachability energy describes the minimum $L^2([0,\infty))$-norm of the control input required to steer the system~\eqref{eq: LTISys} to $\mathbf{x}$ from the origin over an infinite time horizon, whereas the observability energy describes the $L^2([0,\infty))$-norm of the output signal if the system evolves unforced from initial condition $\mathbf{x}$. Canonical BT therefore has the interpretation of preserving the $\mathbf{u}\mapsto\mathbf{y}$ input-output map by projecting the system state onto the subspace of state directions which are simultaneously easiest-to-reach and easiest-to-observe. The canonical balanced truncation reduced model~\eqref{eq: RedLTISys} satisfies the following error bound~\cite[Theorem 7.9]{Antoulas2005Book}:
\begin{align*} 
    \|\mathbf{y}-\mathbf{y}_r \|_{L_2} \leq 2\sum_{k=r+1}^d \sigma_k \, \|\mathbf{u}\|_{L_2},
\end{align*}
where $\sigma_1,\sigma_2,\ldots,\sigma_d$, are the square roots of the eigenvalues of~\eqref{eq: BT-EVP} and are called the \textit{Hankel singular values}~\cite{Antoulas2005Book}.

\subsubsection{Likelihood-informed subspace balanced truncation}\label{subsec: LIS-BT}
The work~\cite{Qian2021Balancing} proposes a balanced truncation model reduction approach applicable to Bayesian inference of the initial state of a linear dynamical system~\eqref{eq: BayesianSystem} as described in \Cref{subsec: problem formulation}. The approach has connections both to canonical control-theoretic balancing (\Cref{subsubsec: control-BT}) and to the likelihood-informed dimension reduction approach (\Cref{subsec: Spantini}), and is based on the following choices of $\mathbf{P}$ and $\mathbf{Q}$:
\begin{align}\label{eq: inference Gramians}
    \mathbf{P} = \mathbf{\Gamma}_{\rm pr}, \qquad \mathbf{Q} = \mathbf{Q}_\epsilon:=\int_0^\infty e^{\mathbf{A}^\top t}\mathbf{C}^\top\boldsymbol{\Gamma}_\epsilon^{-1}\mathbf{C}e^{\mathbf{A}t}\,\text{d}t.
\end{align}
The resulting BT reduced model can be interpreted as preserving state directions which are simultaneously subject to high prior uncertainty while contributing most to the system output signal in the $\mathbf{\Gamma}_{\epsilon}^{-1}$-weighted norm. The reduced model has the form~\eqref{eq: RedLTISys} (with zero input port matrix $\mathbf{B}_r$ in the Bayesian smoothing setting) and defines the following approximate forward operator:
\begin{align}\label{eq: forward_LIS}
    \mathbf{G}_{\rm LIS} = \begin{bmatrix}
        \mathbf{C}_r e^{\mathbf{A}_r t_1} \\
        \vdots \\
        \mathbf{C}_r e^{\mathbf{A}_r t_n }
    \end{bmatrix}\mathbf{V}_r^\top,
\end{align}
which in turn induces the following posterior mean and covariance approximations:
\begin{align}\label{eq: LIS pos both}
\begin{split}
   \PosCov^{\rm LIS} &= (\PrCov^{-1} + \mathbf{G}_{\rm LIS}^\top\ObsCov^{-1}\mathbf{G}_{\rm LIS})^{-1} =  \PrCov - \PrCov \mathbf{G}_{\rm LIS}^\top(\ObsCov+\mathbf{G}_{\rm LIS}\PrCov \mathbf{G}_{\rm LIS}^\top)^{-1}\mathbf{G}_{\rm LIS}\PrCov, \\
   \pmb \mu_{\mathrm{pos}}^{\rm LIS} &=\PosCov^{\rm LIS}\mathbf{G}_{\rm LIS}^\top\ObsCov^{-1}\mathbf m.    
\end{split}
\end{align}
These approximations are generally (tautologically) \textit{sub}-optimal relative to the optimal approximations~\eqref{eq: Spantini PosBoth} defined in~\cite{Spantini2015Optimal}. However, in the context of the Bayesian smoothing problem for a high-dimensional dynamical system, these approximations can lead to greater computational savings than the approach of~\eqref{eq: Spantini PosBoth} because evaluating~\eqref{eq: LIS pos both} only requires evolving the low-dimensional reduced system~\eqref{eq: RedLTISys} whereas evaluating~\eqref{eq: Spantini PosBoth} requires evolving the original high-dimensional system. 

Note that if the `compatibility' condition $\mathbf{A}\mathbf{\Gamma}_{\rm pr} + \mathbf{\Gamma}_{\rm pr} \mathbf{A}^\top \prec \boldsymbol{0}$ holds, the prior covariance $\mathbf{\Gamma}_{\rm pr}$ can be interpreted as a reachability Gramian of the form~\eqref{eq: Gramians} for some choice of input port matrix $\mathbf{B}$. 
Additionally, note that for the Bayesian smoothing problem defined by~\eqref{eq: BayesianSystem}, the matrix $\mathbf{G}^\top\mathbf{\Gamma}_{\rm obs}^{-1}\mathbf{G}$ has the form,
\begin{align}\label{eq: smoothing_Fisher}
    \mathbf{G}^\top\mathbf{\Gamma}_{\rm obs}^{-1}\mathbf{G} = \sum_{i = 1}^n e^{\mathbf{A}^\top t_i}\mathbf{C}^\top\boldsymbol{\Gamma}_\epsilon^{-1}\mathbf{C}e^{\mathbf{A}t_i},
\end{align}
so that in certain continuous-time limits of the discrete observation model~\eqref{eq: BayesianSystem_noise}, the quantity $\mathbf{G}^\top\mathbf{\Gamma}_{\rm obs}^{-1}\mathbf{G}$ converges (up to a scalar multiplying constant) to $\mathbf{Q}_\epsilon$ defined in~\eqref{eq: inference Gramians}. The work~\cite{Qian2021Balancing} exploits these connections between control-theoretic balanced truncation and likelihood-informed dimension reduction to show that the approximations~\eqref{eq: LIS pos both} yield near-optimal posterior covariance approximations for compatible prior covariances when $\mathbf{G}^\top\mathbf{\Gamma}_{\rm obs}^{-1}\mathbf{G}$ approaches $\mathbf{Q}_\epsilon$ as defined in~\eqref{eq: inference Gramians}.

\section{Dimension and model reduction approaches for Bayesian inverse problems with rank-deficient prior covariances}\label{sec: whole_sec_rankdef}
Both the LIS dimension reduction~\cite{Spantini2015Optimal} and the LIS model reduction~\cite{Qian2021Balancing} approaches discussed in \Cref{sec: Background} assume $\PrCov$ is invertible. In this section, we consider linear Bayesian inverse problems with \textit{rank-deficient} prior covariances. \cref{sec: existing methods rankdef} provides new formulations and analyses of the previously discussed LIS dimension and model reduction approaches for rank-deficient prior covariances. \Cref{sec: PD-BT} then presents a new inference-oriented model reduction approach applicable to the Bayesian smoothing problem that naturally handles rank-deficient prior covariances. 

\subsection{LIS dimension and model reduction for rank-deficient prior covariances}\label{sec: existing methods rankdef}
This section provides new LIS dimension and model reduction formulations and analyses applicable to rank-deficient prior covariances. We begin
in \Cref{subsec: 2 Versions gen EVP} by discussing how the LIS projector may be naturally defined for rank-deficient prior covariances via square-root balancing, which is associated with nonsymmetric standard eigenvalue problems instead of symmetric generalized eigenvalue problems. \Cref{subsec: Spantini-rankdef,sec: Qian_rankdef} then discuss theoretical guarantees for LIS dimension and model reduction approaches, respectively, in the case of rank-deficient prior covariances.

\subsubsection{Defining the LIS projector for rank-deficient prior covariances}\label{subsec: 2 Versions gen EVP}
The LIS dimension and model reduction approaches are based on an oblique projector $\mathbf{\Pi}_r$ that is defined in the literature~\cite{Spantini2015Optimal,Qian2021Balancing} by the generalized eigenvalue problem~\eqref{eq: BT-EVP} of the $d$-dimensional matrix pencil $(\mathbf{Q},\mathbf{P}^{-1})$, where  $\mathbf{P} =\mathbf{\Gamma}_{\rm pr}$, and $\mathbf{Q}$ is given by $\mathbf{Q}=\mathbf{G}^\top\mathbf{\Gamma}_{\rm obs}\mathbf{G}$ for LIS dimension reduction and by~\eqref{eq: inference Gramians} for LIS model reduction.

We now consider the case where $\mathbf{P}=\mathbf{\Gamma}_{\rm pr}$ is rank-deficient: let $\textsf{rank}(\mathbf{P})=s$ with $s<d$. 
For simplicity of exposition, we will continue to assume $\mathbf{Q}$ is full-rank, but note that rank-deficient Hessians can be treated similarly as in~\cite{flathFastAlgorithmsBayesian2011}. 
We now discuss how the LIS projector definition and interpretation may be extended to this case of a rank-deficient prior covariance.

We begin by noting that the square-root balancing procedure (cf.~\cref{rem: square root balancing}) provides a natural approach to computing the LIS bases $\mathbf{W}_r,\mathbf{V}_r$ when $\mathbf{P}$ is rank-deficient.
The resulting basis vectors in $\mathbf{W}_r$ are the leading $r$ eigenvectors of
the following nonsymmetric eigenvalue problem:
\begin{align}\label{eq: right_EVP_PQ}
    \mathbf{P}\mathbf{Q}\mathbf{W}=\mathbf{W}\mathbf{\Sigma}^2,
\end{align}
which is equivalent to the generalized eigenvalue problem~\eqref{eq: BT-EVP} when $\mathbf{P}$ is invertible. 
Similarly, the basis vectors in $\mathbf{V}_r$ resulting from the square-root balancing procedure are the $r$ leading eigenvectors of the
associated \textit{left} eigenvalue problem,
\begin{align}\label{eq: left_EVP_PQ}
    \mathbf{V}^\top\mathbf{P}\mathbf{Q}=\mathbf{\Sigma}^2\mathbf{V}^\top, \quad \text{or equivalently}\quad \mathbf{Q}\mathbf{P}\mathbf{V}=\mathbf{V}\mathbf{\Sigma}^2.
\end{align}
Both~\eqref{eq: right_EVP_PQ} and~\eqref{eq: left_EVP_PQ} are well-defined even when $\mathbf{P}$ (or $\mathbf{Q}$) is rank-deficient. 
Observe that the last $d-s$ eigenvalues of~\cref{eq: right_EVP_PQ,eq: left_EVP_PQ} are zero and the square-root balancing procedure yields the normalization $\mathbf{V}_s^\top{\mathbf{P}}\mathbf{V}_s=\mathbf{\Sigma}_s=\mathbf{W}_s^\top \mathbf Q\mathbf{W}_s$, where $\mathbf{\Sigma}_s\in\R^{s\times s}$ is the diagonal matrix whose diagonal entries are the square roots of the nonzero eigenvalues of~\eqref{eq: right_EVP_PQ}, \eqref{eq: left_EVP_PQ} and $\mathbf{W}_s$, $\mathbf{V}_s$ are the associated eigenvectors, respectively.

\subsubsection{Optimality of LIS dimension reduction for linear Bayesian inverse problems with rank-deficient prior covariances}\label{subsec: Spantini-rankdef}
In this section we consider the LIS dimension reduction approach described in \Cref{subsec: Spantini} but with rank-deficient prior covariances. 
We note that the optimality guarantees for the posterior approximations~\eqref{eq: Spantini PosBoth} presented in \Cref{subsec: Spantini} rely on metrics that are defined only for full-rank posterior covariance matrices. However, since the range of the posterior covariance is contained in the range of the prior covariance (one can see this easily from the expression~\eqref{eq: Posterior}), the posterior covariance for a rank-deficient prior covariance will also be rank-deficient. In the case of a rank-deficient prior covariance, we will show that the LIS dimension reduction approach based on the projector defined in \Cref{subsec: 2 Versions gen EVP} inherits optimality guarantees in a lower-dimensional space associated with the range of the prior covariance.

For $\mathbf{P} = \mathbf{\Gamma}_{\rm pr}$ and $\mathbf{Q} = \mathbf{G}^\top\mathbf{\Gamma}_{\rm obs}\mathbf{G}$, let the left and right bases $\mathbf{W}_r,\mathbf{V}_r$  for $r\leq s=\textsf{rank}(\mathbf{\Gamma}_{\rm pr})$ be defined by~\cref{eq: right_EVP_PQ,eq: left_EVP_PQ} as described in \Cref{subsec: 2 Versions gen EVP}. Then, note that $\mathbf{\Pi}_s=\mathbf{W}_s\mathbf{V}_s^\top$ is an oblique projector onto $\textsf{Ran}(\mathbf{\Gamma}_{\rm pr})$. We use this projector to define a restricted measurement model in an $s$-dimensional space associated with the range of the prior covariance, as follows:
\begin{align}\label{eq: reduced obs model 2}
    \mathbf{m} =\mathbf{G}\mathbf{\Pi}_s\mathbf{p}+\pmb{\epsilon} = \mathbf G \mathbf{W}_s \mathbf{V}_s^\top\mathbf{p} + \pmb{\epsilon}= \hat{\mathbf G} \hat{\mathbf{p}} + \pmb{\epsilon},\enspace \qquad  \pmb{\epsilon}\sim \mathcal{N}(\mathbf{0},\ObsCov),
\end{align}
where we have defined $\hat{\mathbf G}:=\mathbf{GW}_s \in \mathbb{R}^{d_{\rm obs}\times s}$, and $\hat{\mathbf{p}}=\mathbf{V}_s^\top \mathbf{p}\in\R^s$ is a restricted unknown state. The prior distribution on $\mathbf{p}$ induces the following prior distribution on the restricted state: $\hat{\mathbf{p}}\sim \cN(\mathbf{0},\hat{\mathbf \Gamma}_{\mathrm{pr}})$, where the restricted prior covariance matrix $\hat{\mathbf \Gamma}_{\mathrm{pr}}:=\mathbf{V}_s^\top{\PrCov}\mathbf{V}_s=\mathbf{\Sigma}_s
$ is now diagonal and invertible. 
The posterior for the problem~\eqref{eq: reduced obs model 2} in the restricted space is then $\hat{\mathbf p}| \mathbf{m} \sim\cN(\hat{\pmb\mu}_{\mathrm{pos}}, \hat{\mathbf \Gamma}_{\mathrm{pos}})$, where
\begin{align}\label{eq: RestrPost}
\begin{split}
    \hat{\mathbf \Gamma}_{\mathrm{pos}} &= \hat{\mathbf \Gamma}_{\mathrm{pr}} - \hat{\mathbf \Gamma}_{\mathrm{pr}}\hat{\mathbf G}^\top(\hat{\mathbf G} \hat{\mathbf \Gamma}_{\mathrm{pr}}\hat{\mathbf G}^\top + \ObsCov)^{-1}\hat{\mathbf G} \hat{\mathbf \Gamma}_{\mathrm{pr}}= (\hat{\mathbf \Gamma}_{\mathrm{pr}}^{-1} + \hat{\mathbf G}^\top\ObsCov^{-1}\hat{\mathbf G})^{-1} \in \mathbb{R}^{s \times s},\\
    \hat{\pmb\mu}_{\mathrm{pos}} &= \hat{\mathbf \Gamma}_{\mathrm{pos}}\hat{\mathbf G}^\top\ObsCov^{-1}\mathbf{m}\in \mathbb{R}^s.
\end{split}
\end{align}
The $d$-dimensional posterior statistics~\eqref{eq: Posterior} for $\mathbf{p} = \mathbf{W}_s\hat{\mathbf{p}}$ can be reconstructed from these $s$-dimensional restricted posterior statistics~\eqref{eq: RestrPost} as follows: $\PosCov = \mathbf{W}_s\hat{\mathbf \Gamma}_{\mathrm{pos}}\mathbf{W}_s^\top$ and $\pmb\mu_{\mathrm{pos}} = \mathbf{W}_s \hat{\pmb\mu}_{\mathrm{pos}}$.
Note that while the full posterior covariance $\PosCov\in\R^{d\times d}$ has rank $s<d$ and is non-invertible, the restricted posterior covariance $\hat{\mathbf{\Gamma}}_{\rm pos}\in\R^{s\times s}$ is invertible. We can therefore apply the LIS dimension reduction approach and theory from \Cref{subsec: Spantini} directly to the restricted inverse problem~\eqref{eq: reduced obs model 2}. 

\begin{corollary}\label[corollary]{cor: restricted optimality}
    Let $\hat{\mathbf \Gamma}_{\mathrm{pos}}^{\mathrm{OLR}}$ and $\hat {\pmb \mu}_{\mathrm{pos}}^{\mathrm{OLR}}$ be the rank-$r$ OLR posterior approximations for the $s$-dimensional restricted Bayesian inverse problem defined by~\eqref{eq: reduced obs model 2}. These approximations optimally approximate the restricted posterior statistics~\eqref{eq: RestrPost} in the sense that:
    \begin{align*}
        \hat{\mathbf \Gamma}_{\mathrm{pos}}^{\mathrm{OLR}}=\arg\min_{\hat{\mathbf{\Gamma}}\in\widehat{\mathcal{M}}_r}d^2_\mathcal{F}(\hat{\mathbf{\Gamma}}_{\rm pos},\hat{\mathbf{\Gamma}}) && and && \hat{\pmb \mu}_{\mathrm{pos}}^{\mathrm{OLR}} =\arg\min_{\hat{\pmb\mu}\in\hat{\mathcal{A}}_r} \mathbb{E} \|\hat{\pmb\mu}_{\rm pos} - \hat{\pmb\mu} \|_{\hat{\mathbf{\Gamma}}_{\rm pos}^{-1}}^2,
    \end{align*}
    where $\widehat{\mathcal{M}}_r = \{\hat{\mathbf \Gamma}_{\mathrm{pr}} - \hat{\mathbf K}_r\hat{\mathbf K}_r^\top\succ 0: \mathrm{rank}(\hat{\mathbf K}_r)\leq r\}$ and $\hat{\mathcal{A}}_r= \{\hat{\mathbf A}_r \mathbf{m} \in \mathbb{R}^s :\mathrm{rank}(\hat{\mathbf A}_r)\leq r\}$.
\end{corollary}

For Bayesian inverse problems with rank-deficient prior covariances, \Cref{cor: restricted optimality} is an analogue of \Cref{thm: Spantini} proving optimality of the LIS approximations in the $s$-dimensional space associated with the range of the prior covariance matrix in which the F\"orstner metric and the Mahalanobis distance based on the restricted posterior precision matrix are well-defined. We now show that the LIS approximations of the $s$-dimensional restricted problem~\eqref{eq: reduced obs model 2} discussed in \Cref{cor: restricted optimality} are equivalent to the LIS approximations of the original $d$-dimensional problem with the rank-deficient prior covariance, in the following sense:
\begin{theorem}\label{thm: LIS Reconstruction}
    Let $\mathbf{W}_r,\mathbf{V}_r\in\R^{d\times r}$, and $\mathbf{\Pi}_r=\mathbf{W}_r\mathbf{V}_r^\top$ be defined by~\cref{eq: right_EVP_PQ,eq: left_EVP_PQ} where $\mathbf{P}=\mathbf{\Gamma}_{\rm pr}$ and $\mathbf{Q}=\mathbf{G}^\top\mathbf{\Gamma}_{\rm obs}^{-1}\mathbf{G}$. Then, let $\mathbf{G}_{\mathrm{OLR}}=\mathbf{G}\mathbf{\Pi}_r$ and let ${\mathbf \Gamma}_{\mathrm{pos}}^{\mathrm{OLR}}\in\R^{d\times d}$ and ${\pmb{\mu}}_{\mathrm{pos}}^{\mathrm{OLR}}\in\R^d$ be the OLR approximations defined in~\eqref{eq: Spantini PosBoth}. These $d$-dimensional posterior approximations are equivalent to the optimal $s$-dimensional approximations $\hat{\mathbf \Gamma}_{\mathrm{pos}}^{\mathrm{OLR}}\in\R^{s\times s}$ and $\hat {\pmb \mu}_{\mathrm{pos}}^{\mathrm{OLR}}\in\R^s$ from \Cref{cor: restricted optimality} in the sense that:
\begin{align*}
    \mathbf W_s \hat{\mathbf \Gamma}_{\mathrm{pos}}^{\mathrm{OLR}} \mathbf W_s^\top={\mathbf \Gamma}_{\mathrm{pos}}^{\mathrm{OLR}}\enspace &\text{and} \enspace \mathbf V_s^\top {\mathbf \Gamma}_{\mathrm{pos}}^{\mathrm{OLR}}\mathbf V_s=\hat{\mathbf \Gamma}_{\mathrm{pos}}^{\mathrm{OLR}}\enspace\text{as well as}\\
    \mathbf W_s \hat{\pmb{\mu}}_{\mathrm{pos}}^{\mathrm{OLR}} ={\pmb{\mu}}_{\mathrm{pos}}^{\mathrm{OLR}} \enspace &\text{and} \enspace \mathbf V_s^\top {\pmb{\mu}}_{\mathrm{pos}}^{\mathrm{OLR}} =\hat{\pmb{\mu}}_{\mathrm{pos}}^{\mathrm{OLR}}.
\end{align*}
\begin{proof}
The proof follows by calculation from the definitions of the OLR approximations of the original and the restricted problem, noting that the restricted prior covariance and Fisher information matrices are equal and diagonal: $\hat{\mathbf \Gamma}_{\mathrm{pr}}=\mathbf{V}_s^\top{\PrCov}\mathbf{V}_s=\mathbf{\Sigma}_s=\mathbf{W}_s^\top \mathbf G^\top\ObsCov^{-1}\mathbf G\mathbf{W}_s =\hat{\mathbf G}^\top\ObsCov^{-1}\hat{\mathbf G}$, due to the normalization of $\mathbf{V}_s$ and $\mathbf{W}_s$.
\end{proof}
\end{theorem}

\Cref{thm: LIS Reconstruction} shows the equivalence between the OLR theory for the original $d$-dimensional inverse problem and for the $s$-dimensional restricted problem~\eqref{eq: reduced obs model 2} defined by using the \textit{oblique} projector $\mathbf{\Pi}_r$ to project the unknown onto $\textsf{Ran}(\mathbf{\Gamma}_{\rm pr})$.
We note that one might naturally consider instead using an \textit{orthogonal} projector to project the unknown onto $\textsf{Ran}(\mathbf{\Gamma})$ and then computing OLR approximations for the resulting restricted problem. However, the resulting approximations would differ from the OLR approximations of the original system and would not enjoy the optimality interpretation of~\Cref{cor: restricted optimality}. The oblique projector defined by square-root balancing is in this sense the more natural way to approach OLR approximation for rank-deficient prior covariances.

\subsubsection{Near-optimality of LIS balanced truncation for Bayesian smoothing problems with rank-deficient prior covariances}\label{sec: Qian_rankdef}
We now consider the LIS balanced truncation approach described in \Cref{subsec: LIS-BT} but with rank-deficient prior covariances. Let $\mathbf{P}$ and $\mathbf{Q}$ be given by~\eqref{eq: inference Gramians}, with $\textsf{rank}(\mathbf{\Gamma}_{\rm pr})=s <d$. Then, for $r\leq s$, let the left and right bases $\mathbf{W}_r,\mathbf{V}_r\in\R^{d\times r}$ be defined by~\cref{eq: right_EVP_PQ,eq: left_EVP_PQ} as described in \Cref{subsec: 2 Versions gen EVP}. These reduced bases then define a reduced model given by~\eqref{eq: RedLTISys} (with $\mathbf{B}_r=\boldsymbol{0}$) which defines an approximate forward operator given by~\eqref{eq: forward_LIS} and posterior approximations given by~\eqref{eq: LIS pos both}. 

The reduced model defined in this way for rank-deficient prior covariances inherits near-optimality guarantees similar to how the rank-deficient LIS dimension reduction approach inherits optimality guarantees as discussed in the previous section. That is, let the restricted $s$-dimensional inverse problem be defined in~\eqref{eq: reduced obs model 2} with solution~\eqref{eq: RestrPost} as before. Denote by $\hat{\mathbf{\Gamma}}_{\rm pos}^{\rm LIS}=\mathbf{V}_s^\top\mathbf{\Gamma}_{\rm pos}^{\rm LIS}\mathbf{V}_s\in\R^{s\times s}$ the posterior covariance in the restricted $s$-dimensional space associated with the LIS balanced truncation covariance approximation $\mathbf{\Gamma}_{\rm pos}^{\rm LIS}$. The reduced model approximation $\hat{\mathbf{\Gamma}}_{\rm pos}^{\rm LIS}$ is again tautologically sub-optimal (in the sense of \Cref{cor: restricted optimality}), but may be cheaper to compute than the optimal approximation in settings where the forward model requires the evolution of the dynamical system. Furthermore, note that rank-deficient prior covariances may still satisfy the prior compatibility condition that $\mathbf{A}\mathbf{\Gamma}_{\rm pr} + \mathbf{\Gamma}_{\rm pr}\mathbf{A}^\top\prec \boldsymbol{0}$. If this condition holds, then the restricted $s$-dimensional posterior covariance associated with the LIS balanced truncation approximation approaches the optimal restricted $s$-dimensional posterior covariance approximation as $\mathbf{G}^\top\mathbf{\Gamma}_{\rm obs}^{-1}\mathbf{G}$ approaches $\mathbf{Q}_\epsilon$ as defined in~\eqref{eq: inference Gramians}. That is, if $\mathbf{\Gamma}_{\rm pr}$ is compatible, then $\hat{\mathbf{\Gamma}}_{\rm pos}^{\rm LIS}\to\hat{\mathbf{\Gamma}}_{\rm pos}^{\mathrm{OLR}}$ in the same continuous-time observation limits considered in~\cite{Qian2021Balancing}.

\subsection{Prior-driven balancing for the linear Bayesian smoothing problem}\label{sec: PD-BT}
We now present \textit{prior-driven balancing}, a new inference-oriented model reduction method applicable to the linear Bayesian smoothing problem introduced in~\Cref{subsec: problem formulation}. The new method is based on interpreting the unknown initial condition as an unknown impulse input to a linear time-invariant control system, and is introduced in \Cref{subsec: PD-BT Intro} and analyzed in \Cref{subsec: PD-BT ApproxQual}.

\subsubsection{Prior-driven balanced truncation}\label{subsec: PD-BT Intro}
We consider the inference of the unknown initial condition $\mathbf{x}(0)=\mathbf{p}$ of the linear dynamical system~\eqref{eq: BayesianSystem_nonoise} from noisy linear measurements~\eqref{eq: BayesianSystem_noise} taken after the initial time. 
Recall that the prior distribution for the unknown initial condition is $\mathbf{p}\sim \mathcal{N}(\mathbf{0},\PrCov)$. 
Let $s=\textsf{rank}(\mathbf{\Gamma}_{\rm pr})$ as before, and let $\mathbf{L}_\mathrm{pr} \in \mathbb{R}^{d \times s}$ denote a square root factor of $\mathbf{\Gamma}_{\rm pr}$ such that $\mathbf{\Gamma}_\mathrm{pr} = \mathbf{L}_\mathrm{pr}\mathbf{L}_\mathrm{pr}^\top$. 
Some algorithms for inference naturally track such a square root factor instead of the prior covariance matrix itself~\cite{LIVINGS2008UnbiasedEnKF,TUANPHAM1998SEEK,Rozier2007ROKF}; the square root factor may alternatively be efficiently computed using iterative methods (see, e.g.,\@ \cite[Chap.\@ 6]{higham:FM2008}).

Our method builds on previous works~\cite{BEATTIE2017inhom,Heink2011inhom,SchroederVoigt23} that have developed system-theoretic balanced truncation methods for inhomogeneous initial conditions. It is also similar to balancing methods for variational data assimilation in the discrete LTI system setting \cite{BoessDiss,GreenDiss}. The key insight is that the system response to a non-zero initial condition is equivalent to the response of a homogeneous system to an impulse input with a suitably designed input port matrix. Note that one way to draw samples of $\mathbf{p}$ is to draw samples of $\mathbf{z}\sim\mathcal{N}(\mathbf{0},\mathbf{I}_s)$ and then let $\mathbf{p}=\mathbf{L}_\mathrm{pr}\mathbf{z}$. We therefore consider the following homogeneous LTI system,
\begin{align}\label{eq: InhomBayesianSys}
    \begin{split}
        \dot{\mathbf{x}}(t)  &= \mathbf{A}\mathbf{x}(t)+\mathbf{L}_\mathrm{pr}\mathbf{u}(t), \enspace \mathbf{x}(0) = \mathbf{0},\\
        \mathbf{y}_\epsilon(t) &= \mathbf{\Gamma}_\epsilon^{-\frac{1}{2}}\mathbf{C}\mathbf{x}(t).
    \end{split}
\end{align} 
We call the system~\eqref{eq: InhomBayesianSys} the \textit{prior-driven} LTI system. 
Let $\delta(t)$ denote the Dirac delta function. The impulse response of the prior-driven system to the impulse input $\mathbf{u}(t)=\delta(t)$ is typically denoted $\mathbf{h}(t)$ and given by $\mathbf{h}(t)=\mathbf{\Gamma}_{\mathbf{\epsilon}}^{-\frac{1}{2}}\mathbf{C}e^{\mathbf{A}t}\mathbf{L}_{\mathrm{pr}}$. 
Then, consider the input $\mathbf{u}(t) = \mathbf{z}\delta(t)$. The response of the prior-driven system to this input is  $\mathbf{y}_\epsilon(t) =\mathbf{h}(t)\mathbf{z} = \mathbf{\Gamma}_\epsilon^{-\frac{1}{2}}\mathbf{C}e^{\mathbf{A}t}\mathbf{L}_\mathrm{pr}\mathbf{z} = \mathbf{\Gamma}_\epsilon^{-\frac{1}{2}}\mathbf{C}e^{\mathbf{A}t}\mathbf{p} = \mathbf{\Gamma}_\epsilon^{-\frac{1}{2}}\mathbf{y}(t)$. Thus, accurate approximation of the input-output map of the prior-driven system~\eqref{eq: InhomBayesianSys} enables accurate approximation of the response $\mathbf{y}(t)$ of the unforced system~\eqref{eq: BayesianSystem_nonoise} with initial condition $\mathbf{p}$, weighted by $\mathbf{\Gamma}_\epsilon^{-\frac{1}{2}}$ to account for the Gaussian observation error in the measurements~\eqref{eq: BayesianSystem_noise}. 
This motivates an alternative balanced truncation formulation for the Bayesian smoothing problem which we now describe.

We assume that the system matrix $\mathbf{A}$ is stable with eigenvalues in the open left half-plane. Then, the infinite reachability and observability Gramians for the prior-driven system~\eqref{eq: InhomBayesianSys} are given by $\mathbf{P}_\mathrm{PD}$ with $\mathbf{A}\mathbf{P}_\mathrm{PD} + \mathbf{P}_\mathrm{PD}\mathbf{A}^\top = -\mathbf{\Gamma}_{\mathrm{pr}}=-\mathbf{L}_\mathrm{pr}\mathbf{L}_\mathrm{pr}^\top$, and $\mathbf{Q}_\epsilon$ as defined in~\eqref{eq: inference Gramians}. 
Note that these Gramians may be rank deficient. Thus, the balanced truncation reduced bases $\mathbf{W}_r,\mathbf{V}_r\in\R^{d\times r}$ should be computed as described in \Cref{subsec: 2 Versions gen EVP} for $\mathbf{P} = \mathbf{P}_{\rm PD}$ and $\mathbf{Q} = \mathbf{Q}_\epsilon$. Then, the reduced bases may be used to obtain the reduced prior-driven system,
\begin{align}\label{eq: ReducedInhomBayesianSys}
    \begin{split}
        \dot{\mathbf{x}}_r(t)  &= \mathbf{A}_r\mathbf{x}_r(t)+\mathbf{L}_{\mathrm{pr},r}\mathbf{u}(t), \enspace \quad\mathbf{x}(0) = \mathbf{0},\\
         \mathbf{y}_{\epsilon,r}(t) &= \mathbf{\Gamma}_\epsilon^{-\frac{1}{2}}\mathbf{C}_r\mathbf{x}_r(t),
    \end{split}
    \end{align}
with $\mathbf{x}_r:= \mathbf{V}_r^\top\mathbf{x}$, and reduced operators $\mathbf{A}_r :=  \mathbf{V}_r^\top\mathbf{A}\mathbf{W}_r$, $\mathbf{L}_{\mathrm{pr},r} :=  \mathbf{V}_r^\top \mathbf{L}_\mathrm{pr}$, and $\mathbf{C}_r := \mathbf{C}\mathbf{W}_r$.

Similar to~\eqref{eq: forward_LIS}, define a reduced forward map constructed from the reduced system~\eqref{eq: ReducedInhomBayesianSys}:
\begin{align*}
    \mathbf{G}_\mathrm{PD} = \begin{bmatrix}
        \mathbf{C}_re^{\mathbf{A}_rt_1}\\ \cdots\\ \mathbf{C}_re^{\mathbf{A}_rt_n}
    \end{bmatrix}\mathbf{V}_r^\top,
\end{align*}
inducing the following approximations to the posterior mean and covariance:
\begin{align}\label{eq: PD-BT PosBoth}
\begin{split}
    \mathbf{\Gamma}_{\mathrm{pos}}^{\mathrm{PD}} &= \PrCov - \PrCov \mathbf{G}_\mathrm{PD}^\top(\ObsCov+\mathbf{G}_\mathrm{PD}\PrCov \mathbf{G}_\mathrm{PD}^\top)^{-1}\mathbf{G}_\mathrm{PD}\PrCov,\\
    \text{and}\quad \pmb{\mu}_{\mathrm{pos}}^{\mathrm{PD}} &= \mathbf{\Gamma}_{\mathrm{pos}}^{\mathrm{PD}}\mathbf{G}_\mathrm{PD}^\top\mathbf{\Gamma}_{\mathrm{obs}}^{-1}\mathbf{m}.
\end{split}
\end{align}
The presented method of \textit{prior-driven balancing} for Bayesian inference (PD-BT) exploits the link between the impulse response of the prior-driven system~\eqref{eq: InhomBayesianSys} and the output of the Bayesian inverse problem~\eqref{eq: BayesianSystem} for reduced posterior computation. 
The new method naturally applies to both full-rank and rank-deficient prior covariances through the use of an appropriate square-root factor as the input port matrix. Additionally, it may be appropriate to use a low-rank approximate square root factor for prior covariances exhibiting strong spectral decay. 

\paragraph{Comparison of computational costs} The main computational difference among the presented approaches---OLR dimension reduction, LIS balanced truncation (LIS-BT), and prior-driven balancing--—lies in the distribution of the computational load. Optimal OLR dimension reduction \cite{Spantini2015Optimal}, only requires solving the nonsymmetric eigenvalue problem \eqref{eq: right_EVP_PQ} or \eqref{eq: left_EVP_PQ} to compute the likelihood-informed subspace. This is typically accomplished by computing an SVD (see square-root balancing in \cref{rem: square root balancing}) and requires $\mathcal{O}(d^3)$ flops \cite{Benner2006NumLinA}. Using an approximate or sparse decomposition of dimension $s$ leads to a more efficient computation of the eigen-decomposition in in $\mathcal{O}(s^2d)$ flops \cite[Section 14.1]{Antoulas2005Book}. However, OLR still requires the action of the full $\mathbf{G}$, which in the cases we are interested in requires evolving the high-dimensional system \eqref{eq: BayesianSystem_nonoise} before the projection onto the LIS is applied. This incurs a cost of $\mathcal{O}(d^3)$ for the standard case of an explicitly given dense system matrix $\mathbf{A}$. This cost can be reduced by orders of magnitude for implicitly given sparse or structured system matrices, which allow for inexpensive matrix-vector products \cite{Benner2008MatrixEqs}.

In contrast, PD-BT involves solving two Lyapunov equations to obtain the reachability and observability Gramians. Then, the resulting nonsymmetric eigenvalue problem is solved via square-root balancing, with costs discussed above. Solving a Lyapunov equation usually requires $\mathcal{O}(d^3)$ flops \cite{Antoulas2005Book}. However, more efficient methods exist for sparse or hierarchical system matrices, as well as when approximate matrix decompositions are available \cite{Benner2006NumLinA}. The second balancing method, LIS-BT \cite{Qian2021Balancing}, requires solving one Lyapunov equation and a nonsymmetric eigenvalue problem. It may also require checking for prior compatibility and solving an additional Lyapunov equation to construct a compatible prior. All of the aforementioned computations are necessary for reduced basis computation. Once the reduced system matrices have been constructed, forwarding the reduced dynamical system \eqref{eq: ReducedInhomBayesianSys} or \eqref{eq: RedLTISys} to build the reduced forward map $\mathbf{G}_\mathrm{PD}$ or $\mathbf{G}_\mathrm{LIS}$ incurs a cost of only $\mathcal{O}(r^3)$ instead of $\mathcal{O}(d^3)$ for $r\ll d$. Thus, the cost of forward model simulations for balancing methods is drastically reduced. Thus, the benefits of model reduction are most significant in settings where the forward model must be evaluated many times, for example when multiple instances of the inverse problem must be solved, e.g., for different realizations of observation data.
    
\subsubsection{Analysis of prior-driven balancing}\label{subsec: PD-BT ApproxQual}
We now prove results concerning the stability and error of the prior-driven balanced reduced model. First, note that the PD-BT reduced system inherits the stability of the full system:
\begin{corollary}[{\cite[Theorem 7.9.1]{Antoulas2005Book}}]
    If the full system~\eqref{eq: InhomBayesianSys} is stable, then the reduced prior-driven LTI system~\eqref{eq: ReducedInhomBayesianSys} is stable as well.
\end{corollary}

We now consider errors in the approximation of the output $\mathbf{y}_{\epsilon}(t)$ of the prior-driven system~\eqref{eq: InhomBayesianSys} and show that our approach inherits error bounds from system-theoretic analyses. 

\begin{proposition}
Consider the input $\mathbf{u}(t) = \delta(t) \mathbf{z}$, where $\mathbf{z}$ is a standard multivariate normal random variable. Let $\mathbf{y}_\epsilon(t)$ and $\mathbf{y}_{\epsilon,r}(t)$ denote the responses of the original~\eqref{eq: InhomBayesianSys} and reduced~\eqref{eq: ReducedInhomBayesianSys} prior-driven LTI systems. Then, 
\begin{align*} 
\mathbb{E}_{\mathbf{z}}\Big[\| \mathbf{y}_\epsilon-\mathbf{y}_{\epsilon,r}\|_{L_2}\Big] \leq 2\sqrt{s} \cdot\sum_{k=r+1}^d \sigma_k.
\end{align*} 
\begin{proof}
    The argument follows the proof of the standard balanced truncation error bound (see, e.g., \cite[Theorem 7.9.3 ]{Antoulas2005Book}) and using that in our case, $\mathbf{z}\sim \mathcal{N}(\mathbf{0},\mathbf{I}_s)$ and thus $\mathbb{E}_{\mathbf{z}}\Big[\|\mathbf{z}\|_2^2\Big] =\mathrm{trace}(\mathbf{I}_s) = s$. With that, $\mathbb{E}_{\mathbf{z}}\Big[\|\delta(t)\mathbf{z}\|_{L_2}\Big]\leq\mathbb{E}_{\mathbf{z}}\Big[\|\mathbf{z}\|_2\Big]\|\delta(t)\|_{L_2} \leq \sqrt{s}$.
\end{proof}
\end{proposition}
We now show how well the impulse response of the reduced prior-driven system~\eqref{eq: ReducedInhomBayesianSys} approximates the impulse response of the full system~\eqref{eq: InhomBayesianSys}. This allows us to derive another bound for errors in the approximation of the output $\mathbf{y}_{\epsilon}(t)$ by its reduced approximation.
Following  \cite[Theorem 3.1 and 3.2]{BEATTIE2017inhom} it holds:
    \begin{proposition}\label[proposition]{prop: SplittingErrorBound}
    Consider the full system~\eqref{eq: InhomBayesianSys} and its reduced version~\eqref{eq: ReducedInhomBayesianSys} obtained by PD-BT with $r\leq s$. Let $\mathbf{S}$ be a solution of the Sylvester equation $$\mathbf{A}^\top\mathbf{S} + \mathbf{SA}_r + \mathbf{C}^\top\mathbf{\Gamma}_\epsilon^{-1}\mathbf{C}_r =\mathbf{0},$$ and let $\bar{\mathbf{S}}\in\mathbb{R}^{(d-r)\times r}$ collect the last $d-r$ rows of $\mathbf{S}$. Similarly, let $\bar{\mathbf{V}},\bar{\mathbf{W}}\in\mathbb{R}^{d\times (d-r)}$ collect the last $d-r$ columns of $\mathbf{V},\mathbf{W}$ and let $\bar{\mathbf{\Sigma}}:=\mathrm{diag}(\sigma_{r+1},\ldots,\sigma_d)\in\mathbb{R}^{(d-r)\times (d-r)}$ be the diagonal matrix of truncated Hankel singular values of $\mathbf{P}_\mathrm{PD}\mathbf{Q}_\epsilon$. Let $\bar{\mathbf{L}}_\mathrm{pr}:=\bar{\mathbf{V}}^\top\mathbf{L}_\mathrm{pr}\in\mathbb{R}^{(d-r)\times s}$ and  $\bar{\mathbf{A}}:=\mathbf{V}_r^\top\mathbf{A}\bar{\mathbf{W}}\in\mathbb{R}^{r\times (d-r)}$.
        Then, for the output functions $\mathbf{y}_\epsilon$ of the full system~\eqref{eq: InhomBayesianSys} and $\mathbf{y}_{\epsilon,r}$ of the balanced, rank-$r$ reduced system~\eqref{eq: ReducedInhomBayesianSys}:
        \begin{align}\label{eq: InhomBayesianErrorBound}
            \mathbb{E}_{\mathbf{z}}\Big[\| \mathbf{y}_\epsilon-\mathbf{y}_{\epsilon,r}\|^2_{L_2}\Big] \leq s \cdot\mathrm{trace}\big[ (\bar{\mathbf{L}}_\mathrm{pr}\bar{\mathbf{L}}_\mathrm{pr}^\top+2\bar{\mathbf{S}}\bar{\mathbf{A}})\bar{\mathbf{\Sigma}}\big].
        \end{align}
\begin{proof}
Denote by
\begin{align*}
    \mathbf{h}(t) = \mathbf{\Gamma}_\epsilon^{-\frac{1}{2}}\mathbf{C}e^{\mathbf{A}t}\mathbf{L}_{\mathrm{pr}}\quad \text{and}\quad \mathbf{h}_r(t) = \mathbf{\Gamma}_\epsilon^{-\frac{1}{2}}\mathbf{C}_re^{\mathbf{A}_rt}\mathbf{L}_{\mathrm{pr},r}
\end{align*}
the impulse responses of the full system~\eqref{eq: InhomBayesianSys} and the balanced, rank-$r$ reduced system~\eqref{eq: ReducedInhomBayesianSys}, respectively. Then, from \cite[Theorem 3.1]{BEATTIE2017inhom} it is known that
\begin{align}\label{eq: InhomImpulseErrorBound}
    \| \mathbf{h}-\mathbf{h}_r\|^2_{L_2} \leq \mathrm{trace}\big[ (\bar{\mathbf{L}}_\mathrm{pr}\bar{\mathbf{L}}_\mathrm{pr}^\top+2\bar{\mathbf{S}}\bar{\mathbf{A}})\bar{\mathbf{\Sigma}}\big].
\end{align}
By construction, it holds that
\begin{align*}
    \| \mathbf{y}_\epsilon-\mathbf{y}_{\epsilon,r}\|^2_{L_2} = \| \mathbf{\Gamma}_\epsilon^{-\frac{1}{2}}\mathbf{C}e^{\mathbf{A}t}\mathbf{L}_{\mathrm{pr}}\mathbf{z}-\mathbf{\Gamma}_\epsilon^{-\frac{1}{2}}\mathbf{C}_re^{\mathbf{A}_rt}\mathbf{L}_{\mathrm{pr},r}\mathbf{z}\|^2_{L_2} \leq  \| \mathbf{h}-\mathbf{h}_r\|^2_{L_2} \|\mathbf{z}\|_2^2.
\end{align*}
Taking the expectation w.r.t. $\mathbf{z}$ yields~\eqref{eq: InhomBayesianErrorBound}.
\end{proof}
\end{proposition}

We now show a result that bounds the difference between the noise-free observations of the full system~\eqref{eq: BayesianSystem} and a reduced version in expectation in an appropriate norm. Let $\mathbf{Y} := \begin{bmatrix}\mathbf{y}(t_1)^\top,\ldots,\mathbf{y}(t_n)^\top\end{bmatrix}^\top \in \mathbb{R}^{d_\mathrm{obs}}$ collect the noise-free outputs of~\eqref{eq: BayesianSystem} and let $\mathbf{Y}_r := \begin{bmatrix}\mathbf{y}_r(t_1)^\top,\ldots,\mathbf{y}_r(t_n)^\top\end{bmatrix}^\top \in \mathbb{R}^{d_\mathrm{obs}}$ with $\mathbf{y}_r(t)=\mathbf{C}_re^{\mathbf{A}_rt}\mathbf{p}_r$ its reduced version. Here, $\mathbf{p}_r=\mathbf{V}_r^\top \mathbf{p} = \mathbf{V}_r^\top\mathbf{L}_\mathrm{pr}\mathbf{z} $, so that $\mathbf{p}_r\sim \mathcal{N}(\mathbf{0},\mathbf{L}_{\mathrm{pr},r}\mathbf{L}_{\mathrm{pr},r}^\top)$. Note that in expectation, the matrices of noise-free measurements $\mathbf{Y}$ and $\mathbf{Y}_r$ can be related to series of impulse responses of the prior-driven system~\eqref{eq: InhomBayesianSys} and its reduced version~\eqref{eq: ReducedInhomBayesianSys}, which motivates to apply \cref{prop: SplittingErrorBound}.

\begin{theorem}\label{thm: output_error_bound}
    Let $\mathrm{rank}(\PrCov) =s$, $\bar{\mathbf{L}}_\mathrm{pr}$, $\bar{\mathbf{S}}$, $\bar{\mathbf{A}}$ and $\bar{\mathbf{\Sigma}}$ as in \cref{prop: SplittingErrorBound} and $\kappa$ a constant, then:
\begin{align}\label{eq: meas_error_bound}
    \mathbb{E}_{\mathbf{p}}\|\mathbf{Y}-\mathbf{Y}_r\|_{\mathbf{\Gamma}_\mathrm{obs}^{-1}}^2 \leq s \cdot \kappa \cdot \mathrm{trace}\big[(\bar{\mathbf{L}}_\mathrm{pr}\bar{\mathbf{L}}_\mathrm{pr}^\top+2\bar{\mathbf{S}}\bar{\mathbf{A}})\bar{\mathbf{\Sigma}} \big].
\end{align}
\begin{proof} We decompose
    \begin{align}
    \mathbb{E}_{\mathbf{p}}\|\mathbf{Y}-\mathbf{Y}_r\|_{\mathbf{\Gamma}_\mathrm{obs}^{-1}}^2 &=\mathbb{E}_{\mathbf{p}}(\mathbf{Y}-\mathbf{Y}_r)^\top\mathbf{\Gamma}_\mathrm{obs}^{-1}(\mathbf{Y}-\mathbf{Y}_r)\notag \\
    &= \mathbb{E}_{\mathbf{p}}\sum_{k=1}^n (\mathbf{y}(t_k)-\mathbf{y}_r(t_k))^\top\mathbf{\Gamma}_\epsilon^{-1} (\mathbf{y}(t_k)-\mathbf{y}_r(t_k))\notag\\
    &= \mathbb{E}_{\mathbf{p}}\sum_{k=1}^n \|\mathbf{\Gamma}_\epsilon^{-\frac{1}{2}}\mathbf{y}(t_k)-\mathbf{\Gamma}_\epsilon^{-\frac{1}{2}}\mathbf{y}_r(t_k)\|_2^2\notag\\
    &=\mathbb{E}_{\mathbf{z}}\sum_{k=1}^n \|\mathbf{\Gamma}_\epsilon^{-\frac{1}{2}}\mathbf{C}e^{\mathbf{A}t_k}\mathbf{L}_{\mathrm{pr}}\mathbf{z}-\mathbf{\Gamma}_\epsilon^{-\frac{1}{2}}\mathbf{C}_re^{\mathbf{A}_rt_k}\mathbf{L}_{\mathrm{pr},r}\mathbf{z}\|_2^2\notag\\
    &\leq \underbrace{\mathbb{E}_{\mathbf{z}}\|\mathbf{z}\|_2^2}_{\substack{=s}}\enspace\sum_{k=1}^n \|\underbrace{\mathbf{\Gamma}_\epsilon^{-\frac{1}{2}}\mathbf{C}e^{\mathbf{A}t_k}\mathbf{L}_{\mathrm{pr}}}_{\substack{\mathbf{h}(t_k)}}-\underbrace{\mathbf{\Gamma}_\epsilon^{-\frac{1}{2}}\mathbf{C}_re^{\mathbf{A}_rt_k}\mathbf{L}_{\mathrm{pr},r}}_{\substack{\mathbf{h}_r(t_k)}}\|_F^2.\label{eq: Sum-Impulse-error}
\end{align}
Note that $\|\mathbf{h}(t_k)-\mathbf{h}_r(t_k)\|_F^2$ is non-negative and finite for all times $t_k$, hence there is a constant $C_k\geq 0$ such that $\|\mathbf{h}(t_k)-\mathbf{h}_r(t_k)\|_F^2 \leq C_k \int_{t_{k-1}}^{t_k} \|\mathbf{h}(t)-\mathbf{h}_r(t)\|_F^2 \mathrm{d}t$ for each $k$ ($t_0 = 0$). Summing up the time steps and setting $\kappa = \max_k C_k$ yields
\begin{align}\label{eq: integral_limit}
    \sum_{k=1}^n\|\mathbf{h}(t_k)-\mathbf{h}_r(t_k)\|_F^2 \leq \kappa \int_{0}^{t_n} \|\mathbf{h}(t)-\mathbf{h}_r(t)\|_F^2 \mathrm{d}t \leq \kappa\| \mathbf{h}-\mathbf{h}_r\|^2_{L_2}.
\end{align}
The right most inequality is justified by the non-negativity of the integrator, s.t.\@ the finite integral is bounded by the infinite integral over $\mathrm{trace}[(\mathbf{h}(t)-\mathbf{h}_r(t))^\top(\mathbf{h}(t)-\mathbf{h}_r(t))]$ yielding the $L_2$-norm of the impulse response. Inserting~\eqref{eq: InhomImpulseErrorBound} and~\eqref{eq: integral_limit} into~\eqref{eq: Sum-Impulse-error} yields the claim~\eqref{eq: meas_error_bound}.
\end{proof}
\end{theorem}
The expected quality of the output approximation of the full Bayesian inference problem by its reduced version---weighted by the known noise precision---thus depends on the truncated Hankel singular values of the prior-driven system~\eqref{eq: InhomBayesianSys}.
\begin{remark}
   The constant $\kappa$ describes the approximation of the integral over the impulse response error by the sum in \eqref{eq: integral_limit}, and can become large if observations are close together, leading to an arbitrarily large output error bound.
   However, we observe in our numerical experiments in \Cref{sec: Numerics} that this bound is quite pessimistic. Thus the primary utility of the result in \Cref{thm: output_error_bound} is to describe the trends in the reduced model error rather than the exact magnitude of the error. This remains true even in the case of closer observation times leading to larger $\kappa$.
\end{remark}

\section{Numerical Experiments}\label{sec: Numerics}
We now present numerical experiments to compare the three dimension and model reduction approaches described in \cref{sec: whole_sec_rankdef} for Bayesian inference of the initial condition of an LTI system. We present the benchmark test problems that we consider in \cref{subsec: test problem} and discuss the results of our experiments in \cref{subsec: Numerics discussion}.

\subsection{Test problems}\label{subsec: test problem}
\paragraph{ISS1R benchmark}
We consider the ISS1R benchmark for system-theoretic model reduction, modeling structural dynamics of the Russian service module of the International Space Station with a discretization of state dimension $d=270$, input dimension $d_\mathrm{in}=3$ and output dimension $d_\mathrm{out}=3$. Its three inputs are the roll, pitch and yaw jets and the three outputs are the roll, pitch and yaw gyroscope measurements \cite{Antoulas2005Book}. The model takes the form~\eqref{eq: LTISys} and the LTI system matrices $\mathbf{A}\in \mathbb{R}^{270 \times 270}$, $\mathbf{B}\in \mathbb{R}^{270 \times 3}$, and $\mathbf{C}\in \mathbb{R}^{3 \times 270}$ are available at \url{http://slicot.org/20-site/126-benchmark-examples-for-model-reduction}. To generate data the system starts at a true initial condition $\mathbf{p}$ drawn from the prior distribution $\mathcal{N}(\mathbf{0},\mathbf{\Gamma}_{\mathrm{pr}})$. We simulate the unforced linear system from time $t = 0$ up to the final measurement time $t_n = 8$, taking measurements at equidistant, discrete times $t_i = 0.1,0.2,\ldots,t_n$. The measurements are generated from the exact evolution of the dynamical system to which we then add $\mathcal{N}(\mathbf{0},\mathbf{\Gamma}_{\epsilon})$-Gaussian measurement noise with $\mathbf{\Gamma}_{\epsilon}=\mathrm{diag}(0.0025^2,0.0005^2,0.0005^2)$. The noise standard deviation is set at around 10\% of the maximum magnitude of the noiseless output.

We provide two experiments using two different types of rank-deficient prior covariances:
\begin{enumerate}
    \item \textit{Incompatible prior covariance}: Empirical covariance of 90 samples drawn from $\mathcal{N}(\mathbf{0},\mathbf{P})$ with $\mathbf{P}$ solving $\mathbf{A}\mathbf{P} + \mathbf{P}\mathbf{A}^\mathrm{T} = -\mathbf{BB}^\mathrm{T}$ and $\mathbf{B}$ given by the ISS1R module benchmark. While $\mathbf{P}$ obtained this way would be a compatible prior, the sample covariance does not satisfy the compatibility condition. This incompatible prior covariance has rank $s = 89$.
    \item \textit{Compatible prior covariance}: We apply the procedure of \cite[Appendix D]{Qian2021Balancing} to modify the above incompatible prior covariance to be compatible. This compatible prior covariance has rank $s = 236$.
\end{enumerate}
We compute the posterior mean and covariance approximations $\mathcal{N}(\pmb{\mu}_{\mathrm{pos}},\mathbf{\Gamma}_{\mathrm{pos}})$ for the initial condition $\mathbf{p}$ using the follow approximation approaches:
\begin{enumerate}
    \item the optimal LIS dimension reduction for rank-deficient prior covariance (OLR, see \cref{subsec: Spantini-rankdef}),
    \item the LIS balanced truncation (LIS-BT) approach described in \cref{sec: Qian_rankdef},
    \item the prior-driven balanced truncation (PD-BT) approach described in \cref{sec: PD-BT}.
\end{enumerate}
We compute 100 independent replicates of the numerical experiment over different realizations of the unknown parameter $\mathbf{p}$ and of the observation noise, and report mean errors over the 100 replicates. 

\paragraph{Two-dimensional heat equation}
As a second example, we consider the heat equation on a two-dimensional plate, and observe the mean temperature in the room. We proceed similarly to the one-dimensional heat equation benchmark example, for which LTI system matrices are available at \url{http://slicot.org/20-site/126-benchmark-examples-for-model-reduction}. This example has been used in previous works \cite{Koe2022TLBTDA,Qian2021Balancing}.

Denote the system matrix for the one-dimensional heat equation by $\mathbf{A}_1\in\mathbb{R}^{50 \times 50}$. Then, the system matrix of the LTI system describing the two-dimensional heat equation with homogeneous Dirichlet boundary conditions is given by $\mathbf{A}=\mathbf{A}_1\otimes\mathbf{I}_{50}+\mathbf{I}_{50}\otimes\mathbf{A}_1\in\mathbb{R}^{2500 \times 2500}$. We let the output matrix $\mathbf{C}\in\mathbb{R}^{1 \times 2500}$ compute the mean temperature in the room.

We consider only an incompatible prior covariance, which is the empirical covariance of 250 samples drawn from $\mathcal{N}(0,\mathbf{R}\mathbf{R}^\top)$ with $\mathbf{R} = \mathrm{diag}(\mathrm{ones}(1,d)) + \mathrm{diag}(0.5\cdot \mathrm{randn}(1,d-1),1) + \mathrm{diag}(0.25\cdot \mathrm{randn}(1,d-2),2)$ and has rank $s=249$. The unforced linear system is simulated from time $t = 0$ up to the final measurement time $t_n = 5$, with measurements at equidistant, discrete times $t_i = 0.2,0.4,\ldots,t_n$ and $\mathcal{N}(0,0.08)$-Gaussian measurement noise added to each observation.

Unless otherwise stated, we proceed as described above for the ISS1R benchmark.

The experiments reported here have been executed on a Dell Latitude 5420 equipped with 32 GB RAM and an 11th Gen Intel(R) Core i7-1185G7 processor. Computations were done in MATLAB 26.1.0.3203278 (R2026a) with Control System Toolbox Version 26.1. running on Ubuntu 24.04.4 LTS Version 1.51.0.

\subsection{Results and discussion}\label{subsec: Numerics discussion}
\begin{figure}[ht]
\centering
\includegraphics[width=0.97\textwidth]{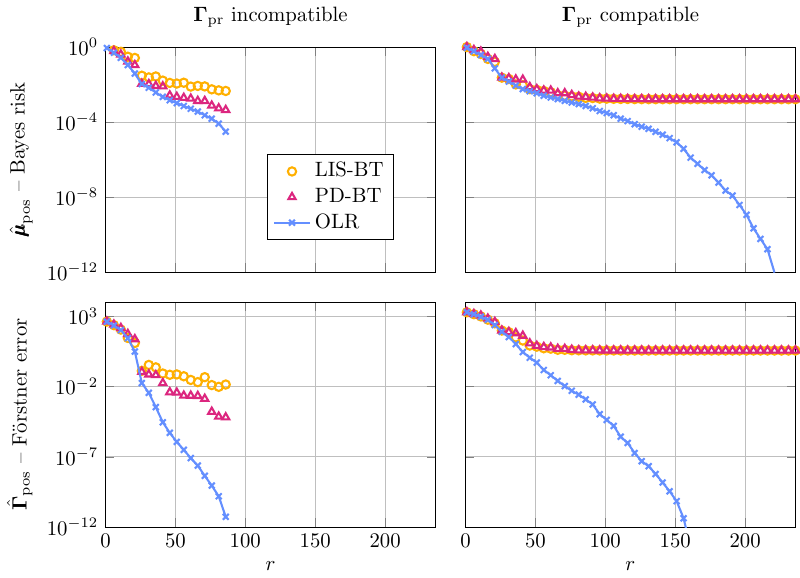}
\caption{Comparison of PD-BT, LIS-BT and optimal dimension reduction (OLR) for the ISS1R module with posterior mean and covariance approximation error in the restricted space, for both an incompatible and a compatible rank-deficient prior covariance.}\label{fig: ISS_reduced_space}
\end{figure}

\paragraph{ISS1R benchmark}
We first compare the posterior approximations in the restricted space associated with $\mathrm{Ran}(\PrCov)$ where LIS dimension reduction is optimal, see \cref{cor: restricted optimality}. The mean and covariance approximation error metrics for both the incompatible and the compatible prior covariance are shown in \cref{fig: ISS_reduced_space}.

The metrics for comparison in the restricted space are the Förstner metric in $\mathbb{R}^s$ for the restricted covariance approximation compared to $\hat{\mathbf{\Gamma}}_{\mathrm{pos}}$ and the $\left\|\cdot\right\|_{\hat{\mathbf{\Gamma}}_{\mathrm{pos}}^{-1}}$-norm Bayes risk for the restricted mean approximation compared to $\hat{\pmb{\mu}}_{\mathrm{pos}}$. To compute these metrics, we restrict to the $s$-dimensional space associated with the range of the prior covariance as discussed in \cref{subsec: Spantini-rankdef}.

For both choices of the prior covariance and for both the mean and covariance approximations, LIS dimension reduction achieves the lowest approximation error. For a compatible prior covariance, there is no difference in approximation quality between LIS-BT and PD-BT for either the mean or the covariance. Both BT methods cannot reach the error of the optimal dimension reduction approximation. With an incompatible prior covariance, however, LIS-BT produces the largest errors. PD-BT can nearly match the optimal dimension reduction approach for the mean approximation with an incompatible prior covariance.
\begin{figure}[ht]
\centering
\includegraphics[width=0.97\textwidth]{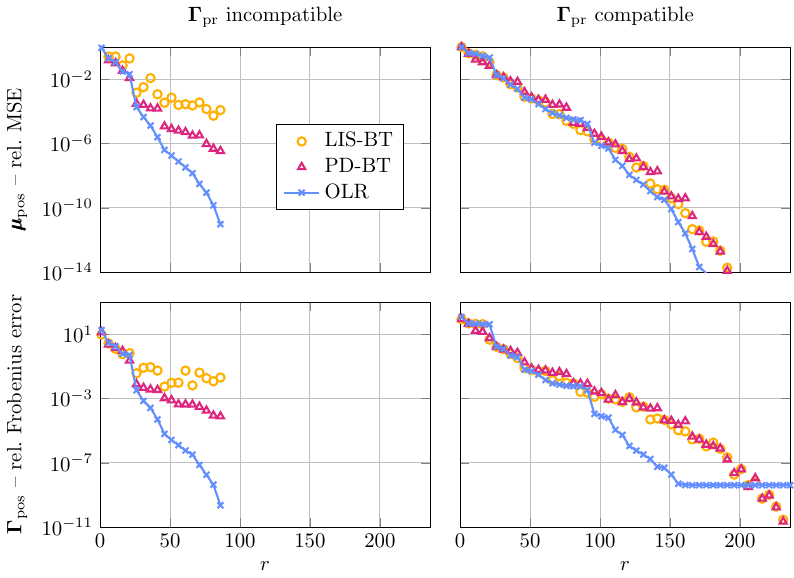}
\caption{Comparison of PD-BT, LIS-BT and optimal dimension reduction (OLR) for the ISS1R module with posterior mean and covariance approximation error in the full space, for both an incompatible and a compatible rank-deficient prior covariance.}\label{fig: ISS_full_space}
\end{figure}

Next, we compare the posterior approximations in the full space $\mathbb{R}^d$. The mean and covariance approximation error metrics for both the incompatible and the compatible prior covariance are shown in \cref{fig: ISS_full_space}.
In the full space, we compare the optimal dimension reduction posterior approximations~\eqref{eq: Spantini PosBoth} to the posterior approximations obtained by PD-BT~\eqref{eq: PD-BT PosBoth} and LIS-BT~\eqref{eq: LIS pos both}. We compute the error to the true posterior covariance $\PosCov$ in the relative Frobenius norm and the relative mean square error (MSE) to the true mean $\pmb{\mu}_\mathrm{pos}$. 
Both balancing approaches behave similarly for a compatible prior. They yield small approximation errors and are close to the LIS approximation for the mean. For an incompatible prior, however, LIS-BT loses its approximation guarantees and yields higher errors than PD-BT. Although PD-BT cannot achieve the error rates of LIS dimension reduction for higher ranks, the relative error remains at $10^{-6}$ and $10^{-4}$ for the mean and covariance approximations, respectively.

\begin{figure}[ht]
\centering
\includegraphics[width=0.97\textwidth]{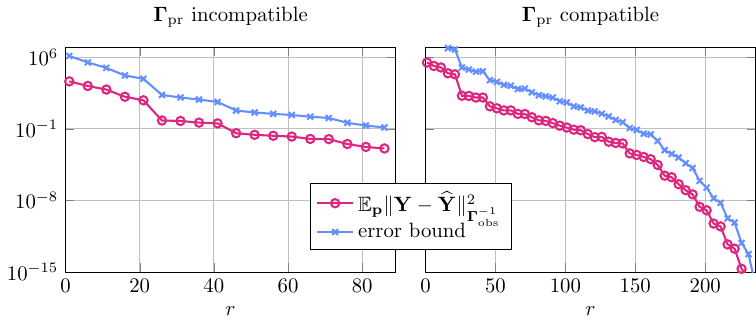}
\caption{Comparison between the output error and the error bound~\eqref{eq: meas_error_bound} for PD-BT given in \cref{thm: output_error_bound}, for both an incompatible and a compatible prior covariance.}\label{fig: ISS_error_bound}
\end{figure}

For the new PD-BT method, \cref{fig: ISS_error_bound} plots the output error bound~\eqref{eq: meas_error_bound} for the true measurements versus the measurements simulated with the PD-BT reduced system.
\Cref{fig: ISS_error_bound} shows that the output error bound~\eqref{eq: meas_error_bound} is tight for both compatible and incompatible priors, with $\kappa\approx10$ in our experimental setting. 

\begin{figure}[ht]
\centering
\includegraphics[width=0.97\textwidth]{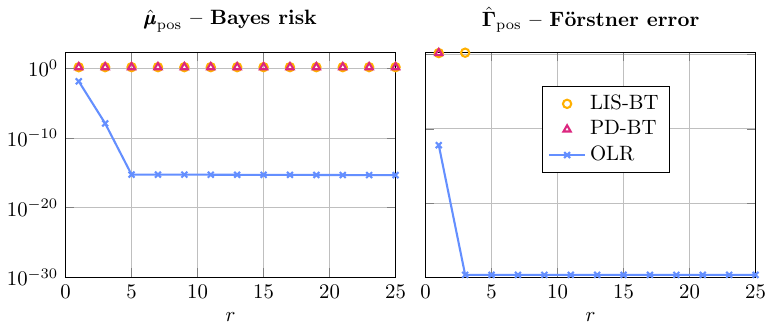}
\caption{Comparison of PD-BT, LIS-BT and optimal dimension reduction (OLR) for the two-dimensional heat equation with posterior mean and covariance approximation error in the restricted space.}\label{fig: heat_posteriors_res}
\end{figure}

\paragraph{Two-dimensional heat equation}
For the second example, we also start by comparing the posterior approximations in the restricted space where LIS dimension reduction is optimal. The mean and covariance approximation error metrics are shown in \cref{fig: heat_posteriors_res}.
The metrics for comparison in the restricted space are again the Förstner metric in $\mathbb{R}^{25}$ for the restricted covariance approximation compared to $\hat{\mathbf{\Gamma}}_{\mathrm{pos}}$ and the $\left\|\cdot\right\|_{\hat{\mathbf{\Gamma}}_{\mathrm{pos}}^{-1}}$-norm Bayes risk for the restricted mean approximation compared to $\hat{\pmb{\mu}}_{\mathrm{pos}}$. To compute these metrics, we restrict to the $25$-dimensional space associated with the leading eigendirections of the eigenvalue problems \eqref{eq: right_EVP_PQ} and \eqref{eq: left_EVP_PQ}, see \cref{subsec: 2 Versions gen EVP}. This restriction is justified because, unlike in the previous example, the forward map $\mathbf{G}\in\mathbb{R}^{25\times 2500}$ is a very wide matrix, since the observations are sparse compared to the high state dimension. Consequently, the Fisher information matrix has a maximum rank of 25, and the OLR approximation based on the eigen-decomposition of $\PrCov\mathbf{G}^\top\ObsCov^{-1}\mathbf{G}$ has at most 25 non-zero eigenvalues. Therefore, we only consider reduced models with $r\leq25$, which is only 1\% of the state dimension.

The heat equation exhibits fast decay of the eigenvalues of the generalized eigenvalue problems for all three approaches. While this results in good approximations of posterior quantities and expected outputs for low ranks, it can also lead to poor posterior covariance approximations and numerical instability. The restricted posterior covariance has a rank of only 5; therefore,  there is no improvement in the OLR approximation after $r=5$ and the approximation errors are practically zero. On the other hand, neither the LIS-BT nor the PD-BT approach yield good approximations of the restricted posterior mean and covariance. This behavior is not surprising given that the LIS-BT posterior approximations only recover the optimality of the OLR posterior approximations for frequent observations. The assumptions for LIS-BT optimality are easily violated, including in the present case with sparse measurements. The second balancing approach, PD-BT, is designed for accurate posterior mean and covariance approximations in the full space with respect to the Frobenius norm, rather than for good posterior approximation in the restricted space. Both balancing methods suffer from numerical instabilities when computing the Förstner metric due to poor conditioning.

\begin{figure}[ht]
\centering
\includegraphics[width=0.97\textwidth]{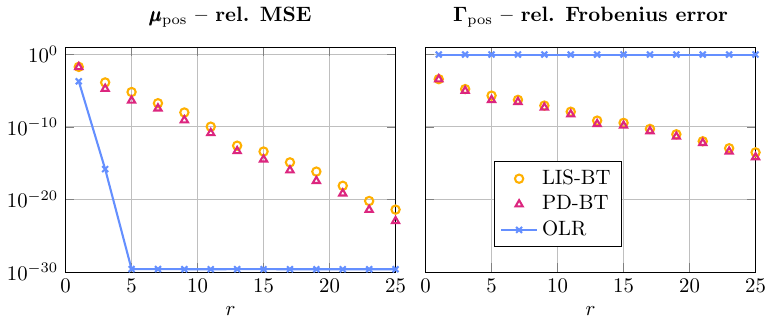}
\caption{Comparison of PD-BT, LIS-BT and optimal dimension reduction (OLR) a for the two-dimensional heat equation with posterior mean and covariance approximation error in the full space.}\label{fig: heat_posteriors}
\end{figure}

Next, we want to compare the errors in the posterior approximations obtained by the three dimension and model reduction methods in the full space $\mathbb{R}^d$, where PD-BT is expected to perform well. We again compute the error to the true posterior covariance $\PosCov$ in the relative Frobenius norm and the relative mean square error (MSE) to the true mean $\pmb{\mu}_\mathrm{pos}$ and plot the results in \cref{fig: heat_posteriors}.

Unlike what we observed in the restricted space, the OLR approach does not accurately approximate the posterior covariance in the full space because it significantly underestimates its rank. Both balancing approaches behave similarly, yielding small approximation errors that decay quickly for the posterior mean and covariance. However, the PD-BT approximation is slightly better. This is to be expected, given the incompatibility of the prior covariance matrix in our experiment. The results presented further demonstrate that PD-BT reliably approximates the posterior mean and covariance in the full space, for which it is designed. Conversely, LIS-BT is designed to recover OLR dimension reduction optimality in the restricted space while avoiding building the full forward model. The presented example has shown that LIS-BT may fail to recover the OLR approximation under ill-conditioning and sparse measurements.

\begin{figure}[ht]
\centering
\includegraphics[width=0.6\textwidth]{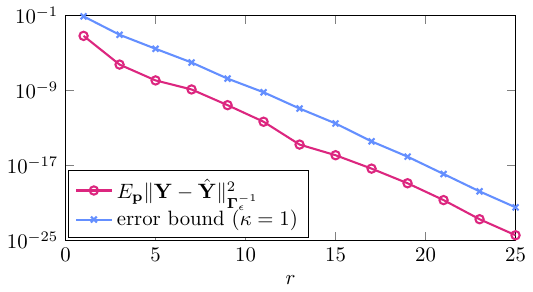}
\caption{Comparison between the output error and the error bound~\eqref{eq: meas_error_bound} with $\kappa = 1$ for PD-BT given in \cref{thm: output_error_bound}, for the two-dimensional heat equation.}\label{fig: heat_error_bound}
\end{figure}

For the new PD-BT method, \cref{fig: heat_error_bound} plots the output error bound~\eqref{eq: meas_error_bound} for the true measurements versus the measurements simulated with the PD-BT reduced system. The matrix exponential of the system matrix $\mathbf{A}$ for the heat equation decays very quickly, causing the bound to hold with $\kappa=1$, as observed in the plot.

The experiments demonstrate that PD-BT is advantageous for incompatible, rank-deficient prior covariances. In these cases, PD-BT yields an accurate and inexpensive forward model and posterior approximation. LIS dimension reduction is optimal in directions with non-zero prior uncertainty but requires computationally expensive full forward models. The use of inexpensive reduced models in both BT approaches alleviates this issue, providing a good approximation in the full space despite being theoretically suboptimal.

\section{Conclusion}\label{sec: Conclusions}
We have proposed three new methods for approximating posterior distributions in linear Gaussian Bayesian inference problems with rank-deficient prior covariances. The first two new methods reformulate a dimension reduction and a balanced truncation model reduction approach based on likelihood-informed subspaces for the case of rank-deficient prior covariances. We show how the optimality guarantees initially developed for full-rank prior covariances in~\cite{Spantini2015Optimal} are applicable to this new setting.

The third new proposed method is \textit{prior-driven balancing}, which is a model reduction method applicable to Bayesian inference of the initial condition of a linear dynamical system from linear measurements taken after the initial time. In contrast to the LIS-based balanced truncation approach~\cite{Qian2021Balancing}, where the prior covariance serves as an analogue of a reachability Gramian, the new approach uses a square root factor of the prior covariance to define an input port matrix for a Gaussian impulse input. We show that the impulse response of this prior-driven LTI system corresponds to the noiseless output signal in the Bayesian inverse problem of interest. This approach inherits error bounds and stability guarantees from standard BT and is applicable for both full-rank and rank-deficient prior covariances. Numerical experiments demonstrate that the proposed PD-BT approach yields accurate and efficient posterior approximations, especially in cases where LIS-BT fails to provide robust and accurate results. 

Several directions are of interest for future work, including (i) analyzing the approximation error of the posterior mean and covariance approximations of the PD-BT approach, (ii) investigating the application of other system-theoretic model reduction methods, including time-limited balanced truncation~\cite{koenig2022TLBT4DVAR,Kuerschner2018TLBT,Wodek1990TLGramians} or interpolatory methods~\cite{AntBeaGuc2020InterpolatoryMeths,Benner2017MORBook,GugercinIRKA2008}, to the prior-driven system, and (iii) extending the approach to nonlinear forward models and non-Gaussian noise and prior distributions.

\backmatter
\section*{Statements and Declarations}
\paragraph*{Funding}
J.K. and M.F. have been partially funded by the Deutsche Forschungsgemeinschaft (DFG) --- Project-ID 318763901 -- SFB1294. E.Q.\ was supported by the United States Air Force Oﬃce of Scientific Research (AFOSR) award FA9550-24-1-0105 (Program Oﬃcer Dr.\ Fariba Fahroo). J.K. gratefully acknowledges financial support by the Fulbright U.S. Student Program, which is sponsored by the U.S. Department of State and the German-American Fulbright Commission. The contents of this paper are solely the responsibility of the authors and do not necessarily represent the official views of the Fulbright Program, the Government of the United States, or the German-American Fulbright Commission.

\paragraph*{Conflict of interest statement}
The authors have no relevant financial or non-financial interests to disclose. E.Q.\@ is an Editorial Board Member for the topical collection \textit{Model Reduction and Surrogate Modeling 2024 (MORe24)}.

\paragraph*{Code availability statement}
All code was implemented in MATLAB and is available at \url{https://github.com/joskoUP/PD-BT}. The plots in the paper were created using the version with the Git commit hash \texttt{f4121a4}, which was updated on May 10, 2026.

\end{document}